\documentclass[11pt,letterpaper]{article}
\usepackage[T1]{fontenc}
\usepackage[utf8]{inputenc}
\usepackage{lmodern}
\usepackage[margin=1in]{geometry}
\usepackage{amsmath,amssymb,amsthm,mathtools}
\usepackage{aliascnt,needspace}
\usepackage{array,booktabs,longtable}
\usepackage{enumitem,microtype,xcolor}
\usepackage{xurl}
\usepackage[colorlinks=true,linkcolor=blue!45!black,citecolor=blue!45!black,urlcolor=blue!45!black]{hyperref}
\usepackage[nameinlink,noabbrev]{cleveref}
\usepackage{listings}
\usepackage{fancyhdr}
\setlist{topsep=4pt,itemsep=3pt,parsep=1pt}

\newtheorem{theorem}{Theorem}[section]
\newaliascnt{lemma}{theorem}
\newtheorem{lemma}[lemma]{Lemma}
\aliascntresetthe{lemma}
\newaliascnt{proposition}{theorem}
\newtheorem{proposition}[proposition]{Proposition}
\aliascntresetthe{proposition}
\newaliascnt{corollary}{theorem}
\newtheorem{corollary}[corollary]{Corollary}
\aliascntresetthe{corollary}
\theoremstyle{definition}
\newaliascnt{definition}{theorem}

\aliascntresetthe{definition}
\newaliascnt{remark}{theorem}

\aliascntresetthe{remark}

\crefname{theorem}{Theorem}{Theorems}
\Crefname{theorem}{Theorem}{Theorems}

\crefname{lemma}{Lemma}{Lemmas}
\Crefname{lemma}{Lemma}{Lemmas}

\crefname{proposition}{Proposition}{Propositions}
\Crefname{proposition}{Proposition}{Propositions}

\crefname{corollary}{Corollary}{Corollaries}
\Crefname{corollary}{Corollary}{Corollaries}

\crefname{definition}{Definition}{Definitions}
\Crefname{definition}{Definition}{Definitions}

\crefname{remark}{Remark}{Remarks}
\Crefname{remark}{Remark}{Remarks}

\crefname{equation}{equation}{equations}
\Crefname{equation}{Equation}{Equations}

\crefname{section}{Section}{Sections}
\Crefname{section}{Section}{Sections}
\crefname{subsection}{Section}{Sections}
\Crefname{subsection}{Section}{Sections}

\crefname{table}{Table}{Tables}
\Crefname{table}{Table}{Tables}

\newcommand{\Z}{\mathbb Z}
\newcommand{\R}{\mathbb R}
\newcommand{\T}{\mathbb T}
\newcommand{\A}{\mathcal A}
\newcommand{\M}{\mathfrak M}
\newcommand{\m}{\mathfrak m}
\newcommand{\e}{\mathrm e}
\newcommand{\dd}{\,\mathrm d}
\newcommand{\eps}{\varepsilon}

\newcommand{\Si}{\mathfrak S}
\newcommand{\Ji}{\mathfrak J}
\DeclareMathOperator{\rad}{rad}
\DeclareMathOperator{\lcm}{lcm}
\title{\textbf{Linear bounds and restricted moments\\ in Waring's problem}}
\author{Morley Davidson}
\date{27 November 2026}
\hypersetup{pdfauthor={Morley Davidson},pdftitle={Linear bounds and restricted moments in Waring's problem}}
\begin{document}
\maketitle
\begin{abstract}
We prove $G(k)<6.873k$ for every positive integer $k$ and
$\limsup_{k\to\infty}G(k)/k\le6.87205868$. For eighth and ninth powers
we obtain $G(8)\le38$ and $G(9)\le44$. We estimate moments of Weyl
sums over integers with marked prime factors by bounding the sum of
squared large values on each fibre of a power map. Vinogradov's mean
value theorem and arithmetic van der Corput transformations control the
resulting Gram matrices. For $k=2^j$, $j\ge3$, we obtain $G_1(k)=4k$
and compute the density of the represented integers for a range of
smaller numbers of summands. At the same linear threshold as for $G(k)$,
the roots can be squarefree with a common number of prime factors.
A reciprocal exponent sum of $6.873$ suffices to represent all large
integers by mixed powers. In particular, the number of consecutive
exponents needed, starting at $k$, has order $k$. For large $k$ we obtain
$\min_{1\le n\le P}\|\alpha n^k\|\le P^{-1/(12.04k)}$ uniformly in
$\alpha\in\mathbb R$. Finally, $\gg A^{0.91709477}$ integers up to $A$
have a growing number of representations by three positive cubes with
distinct squarefree roots. The full weighted Hooley--Manin asymptotic
would give lower density greater than $0.09521$ for sums of three
positive cubes.
\end{abstract}
{\small\setcounter{tocdepth}{1}\tableofcontents}
\section{Introduction}
Let $G(k)$ be the least integer $s$ for which every sufficiently large
positive integer is a sum of at most $s$ positive integral $k$th powers.
Write $G_1(k)$ for the analogous threshold with a set of exceptions of
natural density zero, and $G^+(k)$ for the least count giving positive lower density. Davenport proved $G(4)=16$ and $G_1(3)=4$ \cite{DavenportFour,DavenportThree}. The obstruction to three cubes on the classes $4,5\pmod9$ gives the lower bound in the latter assertion. The variable $k$ is fixed when a representation
estimate is stated. Implicit constants and the point at which the estimate
takes effect may depend on $k$.

Hardy and Littlewood's circle method separates the congruence restrictions
from the analytic cost of representation. Vinogradov's use of smooth
numbers and the efficient differencing of Wooley \cite{W92,W95} reduced
the leading coefficient in bounds of order $k\log k$. Br\"udern and
Wooley \cite{BW23} obtained $G(k)\le\lceil k(\log k+4.20032)\rceil$.
Chow and Dietmann \cite[Theorem 1.1]{CD} proved $G(k)\ll k$. Their
first version announces $G(k)\le8.4k$ after the theorem and gives the
details for the order bound. We use their construction of interlaced prime
intervals and their Gram-matrix argument. We retain the magnitude on each
fibre, and iterate the resulting moment estimate directly. Completion on
successive differences and on their Fourier transforms then supplies an
explicit range of admissible moduli.
The full mean-value theorem of Bourgain, Demeter and Guth \cite{BDG},
and Wooley's proof by nested efficient congruencing \cite{W19}, enter the
degree reduction and the intermediate-power estimates.
Several of our numerical bounds depend on a choice of parameters, and we
have not sought the best choice in every case.

\begin{theorem}\label{thm:main}
For every $k\ge1$ one has $G(k)<6.873k$. Moreover
\[
 \limsup_{k\to\infty}\frac{G(k)}k
 \le \frac{34360293397827451}{5000000000000000}
 =6.8720586795654902 .
\]
The bounds in the third column of \cref{tab:small} hold.
\end{theorem}
\begin{table}[ht]
\centering\small
\begin{tabular}{rrrr|rrrr}
\toprule
$k$&Earlier&$G(k)\le$&$t$&$k$&Earlier&$G(k)\le$&$t$\\
\midrule
8&39&38&20&15&97&83&42\\
9&47&44&23&16&105&90&46\\
10&55&51&26&17&113&97&49\\
11&63&57&29&18&121&104&52\\
12&72&64&33&19&129&110&56\\
13&81&70&36&20&137&117&59\\
14&89&77&39&&&&\\
\bottomrule
\end{tabular}
\caption{Bounds for intermediate powers. The earlier bounds are from
\cite[Theorem 1.1]{W16} and \cite[Theorem 1.3]{BW23}. At order $t$ the
exceptional set on primitive locally soluble classes has a power saving.}
\label{tab:small}
\end{table}
For example, the first row gives two different conclusions about eighth
powers. Every sufficiently large positive integer is a sum of at most
$38$ eighth powers. With exactly $20$ summands one still represents
every integer up to $N$ that has a primitive local representation at
each prime, apart from $O(N^{1-\epsilon})$ exceptions for some
$\epsilon>0$.

\begin{theorem}\label{thm:density}
For each $k\ge21$ there is an even integer $S$ with $S+1<6.873k$.
For $t\ge S/2$, all but $O_{k,t}(N^{1-\epsilon_{k,t}})$ integers
$n\le N$ having a primitive local representation everywhere are sums of
exactly $t$ positive $k$th powers. Here a primitive local representation
has a unit summand at the given prime and $\epsilon_{k,t}>0$.
For $8\le k\le20$ the same assertion holds from the indicated $t$ in
\cref{tab:small} onwards. In particular
\[
 G_1(k)<3.4365k\quad(k\ge21,\ k\text{ not a power of two}),
 \qquad G_1(2^j)=4\cdot2^j\quad(j\ge3).
\]
If $k=2^j$, $j\ge3$, and $t$ is in the indicated range below $4k$,
the number of integers $n\le N$ represented by exactly $t$ powers is
\begin{equation}\label{eq:dyad-density}
 \frac{t}{4k(1-2^{-k})}N+O_{k,t}(N^{1-\epsilon_{k,t}}).
\end{equation}
\end{theorem}
The same sharp moment gives $G^+(k)\le S/2$ for every $k\ge21$,
including powers of two. Indeed, for $t=S/2$ the marked set supplies
$\gg P^t$ tuples and \eqref{eq:sharp} bounds their collision count by
$O(P^{2t-k})$. Cauchy's inequality gives $\gg P^k$ distinct sums in
$[1,tP^k]$.
\begin{theorem}\label{thm:primes}
Let $k\ge21$ and take $s=S+1$ as in \cref{thm:density}. For every
$\eta>0$ there are $c>0$ and an integer $\ell$ such that every sufficiently
large $n$ soluble as a sum of $s$ unit $k$th powers in $\Z_p$ at each
prime $p$ is a sum of $s$ positive $k$th powers whose roots are squarefree
products of exactly $\ell$ primes, all between $n^c$ and $n^\eta$.
The roots may be taken pairwise coprime.
The number of ordered representations of this kind is
$\gg n^{s/k-1}(\log n)^{-s}$. At $t=S/2+1$ the corresponding conclusion
holds outside $O_A(N(\log N)^{-A})$ locally admissible targets $n\le N$
for each $A>0$.

For all sufficiently large $n$, the same lower bound at $s=S+1$
holds with squarefree roots having a common fixed number $\ell^*$ of
distinct prime factors, each factor belonging either to a fixed finite
set or to $[n^c,n^\eta]$.
\end{theorem}
For powers of two, the values $G_1(8)=32$, $G_1(16)=64$ and
$G_1(32)=128$ were obtained by Vaughan and Wooley. Br\"udern and
Wooley gave progression estimates for sixteenth and thirty-second powers
\cite[Theorems 9.1 and 9.2]{BW23}. Formula \eqref{eq:dyad-density}
extends the density calculation to every $2^j$ with $j\ge3$.

The efficient-differencing proofs of $G(8)\le38$ and $G(9)\le46$
were recorded in a working manuscript dated 22 November 2026 \cite{working}.\footnote{It
was evidently a busy Sunday for Waring's problem. The PDF of this working
draft was compiled about 43 minutes before Chow and Dietmann submitted
their preprint to arXiv.}
They are independent of the factors-and-fibres argument used here,
which gives a second proof for eighth powers and lowers the ninth-power
bound to $44$. We sketch the earlier arguments in \cref{sec:small}.

In \cref{sec:bilinear} we use the unrestricted integer at the end of each
marked factorization to pass to almost-prime roots. \Cref{sec:systems}
applies the same sharp moment to several simultaneous equations.
\Cref{sec:fractional} gives a uniform estimate for small fractional parts, and \cref{sec:mixed} treats exponents of different degrees.
The final section considers the much sparser problem of three positive cubes.

The accompanying Lean development proves the earlier representation
bounds $G(8)\le38$ and $G(9)\le46$. For the present bound $G(9)\le44$,
it checks the finite moment deductions from the stated analytic
estimates. The formal deduction of the uniform coefficient $6.873$
also takes the initial moments, shortening estimates and arc bounds as
hypotheses. Their analytic proofs are given here. The numerical
certificates \cite{cert} check the finite inequalities by exact arithmetic.

\paragraph{Acknowledgement.}
Some thirty years ago, as a doctoral student at the University of Michigan,
I had the good fortune to learn about Waring's problem from Robert Vaughan
and my advisor, Trevor Wooley. I remain grateful to them both for their
teaching and guidance, and for the clarity and beauty of their writing
on the subject. In revising this preprint, I hope to bring its arguments
into clearer order and to give readers something of the pleasure I found
in their papers.

\paragraph{Use of generative AI.}
OpenAI's ChatGPT (GPT-6) was used extensively for mathematical exploration,
the development and refinement of proofs, numerical experimentation,
preparation of certificates, assistance with Lean formalization,
literature research, and drafting and revision of the manuscript.
The author assumes responsibility for the mathematical claims, references,
and final contents of this article.
\section{Marked roots and their complete sums}\label{sec:roots}
Put $\e(t)=\exp(2\pi i t)$ and give $\T=\R/\Z$ its probability measure.
We use $P$ for the size of a root and $X$ for its size after some prime
factors have been removed. For $R\ge2$ put
\[
 \A(X,R)=\{1\le n\le X:p\mid n\Longrightarrow p\le R\},\qquad
 f_R(\alpha;X)=\sum_{n\in\A(X,R)}\e(\alpha n^k).
\]
A bound with an admissible excess $E_u$ at order $u$ has the form
$\int|f_R(\alpha;X)|^u\dd\alpha\ll X^{u-k+E_u+\eps}$
for $R\le X^{\eta_0(\eps,k,u)}$. At even orders the same bound holds
on every subset of the roots by orthogonality.
We use only finitely many bounds for any fixed $k$ and final accuracy.

Fix a small $\gamma>0$. Choose distinct small positive exponents $a_i$
($1\le i\le R_0$) of sum $1-2\gamma$ and exponents $b_j$
($1\le j\le J$) of sum $2\gamma-\xi$, where
$0<4\xi<\min b_j\le\max b_j<\min a_i$.
All exponents may be made smaller than any prescribed positive multiple
of $\gamma$. To see this, first choose many $a_i$ close to their average.
Choose many more $b_j$ with average well below $\min a_i$ and finally
adjust one of them by $-\xi$. Choose $\omega>0$ so that the prime
intervals
\[
 I_v=(P^v,P^{(1+\omega)v}]\qquad(v=a_i\ \hbox{or}\ b_j)
\]
are disjoint and $\omega<\xi/4$. Require also
$\log(1+\omega)<1/2$. A full marked product $d$ has one prime from each
interval. The root set consists of $d y$ with $1\le y\le P/d$.
All its coefficients are one. Indeed $P^{\xi/2}<P/d\le P^{2\xi}$, so
$y$ is smaller than every marked prime. The intervals identify the
individual primes in the factorization. To see the point with just two
marks, write $n=p_1p_2y$, with $p_i$ in disjoint prime intervals
$I_i$ and $y$ smaller than every prime in either interval. Factoring
$n$ recovers each $p_i$ from its interval. When $p_1$ is fixed and
removed, the remaining root is $p_2y$, and the phase
$\alpha n^k$ becomes $(p_1^k\alpha)(p_2y)^k$. Thus a fixed value of
$p_1^k\alpha$ collects several original frequencies on one fibre.

Write $F(\alpha)$ for its Weyl
sum. Mertens' theorem, applied separately to the fixed intervals, gives
\begin{equation}\label{eq:mass}
 c_P=\sum_d d^{-1}\longrightarrow (\log(1+\omega))^{R_0+J}>0,
 \qquad \sum_m m^{-1}\le1
\end{equation}
for large $P$, where $m$ ranges over the products from any fixed subset
of intervals. The second inequality includes the empty product.

Only the $a_i$ will be removed in the moment argument. For a set $I$ of
remaining $a_i$ and a fixed choice of primes at the removed indices put
$X=P/(\text{product of removed primes})$. The state $F_I(\alpha;X)$
has one prime at each remaining primary interval, all auxiliary primes,
and the same bottom cofactor bound $X/(\text{product of remaining marks})$.
In particular $P^\gamma\le X\le P$. By taking all window exponents
sufficiently small the roots in every state belong to the smooth set for
all the prescribed seed estimates. The unrestricted cofactor is smaller
than even the smallest marked prime at every state.

For $q\in\mathbb N$ put
\[
 S_k(q,a)=\sum_{r=1}^q\e(ar^k/q),\qquad
 V(\beta)=\int_0^P\e(\beta x^k)\dd x,\qquad
 v(b)=P^{-1}V(bP^{-k}).
\]
Thus $v(b)=\int_0^1\e(bx^k)\dd x\ll_k(1+|b|)^{-1/k}$.
Write $\M(Q)$ for the union of the arcs $q\le Q$,
$(a,q)=1$, $|q\alpha-a|\le QP^{-k}$. We may also use their rectangular
enlargement $\mathcal K(Q)$ with $|\alpha-a/q|\le QP^{-k}$.
For $Q=P^\delta$ with $\delta$ sufficiently small these rectangular
arcs are disjoint.
\begin{lemma}\label{lem:core}
For some $\delta>0$ the estimate
\[
 F(a/q+\beta)=c_P q^{-1}S_k(q,a)V(\beta)+O(P^{1-\xi/2+2\delta})
\]
holds for $q\le P^\delta$ and $|\beta|\le P^{-k+\delta}$.
\end{lemma}
\begin{proof}
All $d$ are units modulo $q$. Multiplication by $d$ permutes the
residues, so the mean of $\e(a(dy)^k/q)$ over $y\bmod q$ is
$S_k(q,a)/q$. Its error in any initial interval is at most $2q$.
Partial summation on $1\le y\le P/d$ and the change $x=dy$ give the
displayed main term for that $d$, with factor $1/d$ and error
$O_k(q(1+P^k|\beta|))$. There are at most
$P^{1-\xi/2}\sum_d1/d$ marked products.
\end{proof}

\section{Conductor estimates}\label{sec:conductor}
Let $m$ be squarefree, with at most $R_0$ prime factors, all
$\le X^\eta$ and tending to infinity. Suppose $q\mid m^k$,
$m\le X^{M/k}$, $1<U=\log_Xq\le M<k$. In this section the sum
\begin{equation}\label{eq:gram-entry}
 \sum_{n\in I,(n,m)=1}\e(an^k/q),\qquad (a,q)=1,
\end{equation}
is taken over any integer interval of length at most $X$. Put
$\theta=M/k$. We need uniform estimates $\ll X^{1-\sigma+o(1)}$,
where the error in the exponent can also contain $O_{k,R_0}(\eta)$.
In the estimates below, rounding a greedy divisor costs $O_k(\eta)$
in the exponent, independently of $R_0$. Implied multiplicative
constants may depend on the number of marks.

\begin{lemma}\label{lem:weyl}
For $2\le j\le k$ an available saving is
\[
 \max\{2^{1-j},1/(j(j-1))\}
 \min\{U/(j+1)-j\theta,\ 1-U/(j+1)-\theta,\ j-U-j\theta\}.
\]
\end{lemma}
\begin{proof}
Split into residue classes modulo $B$, with
$v_p(B)=\max(1,\lceil v_p(q)/(j+1)\rceil)$ at $p\mid m$.
Write $n=c+B t$ on a unit class. All coefficients of degree greater
than $j$ are integral. The coefficient of $t^j$ has denominator
$q_j=\prod_{p\mid q}p^{\max(0,v_p(q)-jv_p(B))}$, since
$p>k$ and $p\nmid c$. In particular
$X^{U/(j+1)-j\theta}\le q_j\le X^{U/(j+1)}$ and the available length
is $Y\asymp X/B$. The usual exact-rational Weyl inequality, in its
classical form or its mean-value form \cite[Theorem 5]{Bourgain}, gives
a normalized bound $\ll X^\eps(q_j^{-1}+Y^{-1}+q_jY^{-j})^\nu$
with $\nu$ the larger of the two displayed exponents. For shorter
intervals the unnormalized bound increases with $Y$
($j\nu\le1$). The remaining factor satisfies
$q_jY^{-j}\ll X^{U-j+j\theta}$, also directly from
$q_j\le q^{1/(j+1)}$ and $B\le q^{1/(j+1)}m$.
Sum over the $B$ classes.
\end{proof}

\begin{lemma}\label{lem:diffcells}
For $0\le j\le k-2$ and $2^j\sigma\le1/2$ the saving
$\sigma-O_k(\eta)-\eps$ is available throughout
\[
 j\theta+2^{j+1}\sigma\le U\le j+2-(2^{j+2}-2)\sigma .
\]
For $j>0$ assume $\theta\le1-2^j\sigma$.
\end{lemma}
\begin{proof}
At conductor primes allow polynomial phases reducing to a nonzero scalar
times $n^v$ modulo $p$, with $2\le v\le k$. At primes of $m$ away
from the conductor allow finitely many excluded classes. Completion
saves $g_0(U)=\max(0,\min(U/2,1-U/2))$, up to $X^\eps$. Indeed each
complete Fourier coefficient on units is $O_k(p^{e/2})$: for $e=1$
use Weil's bound, and for $e\ge2$ the stationary roots modulo $p$
lie in at most $v-1$ unit classes with invertible second derivative.
The prime-power calculation is given below. At primes away from $q$
remove the exclusions by inclusion--exclusion. Substitution along the
resulting progressions changes $n$ affinely with invertible slope at
every conductor prime and gives intervals of length at most $X$.
The resulting completion bound is $O_{k,R_0}((X/\sqrt q+\sqrt q)\log(2q))$.

For $\theta\le b\le\min(1,U)$ take $B\mid q$ divisible by
$\rad(q)$ with $X^b\le B\le X^{b+\eta}$. Start at $\rad(q)$ and
adjoin one available prime factor at a time. Averaging shifts $Bh$ with
$h\le H=\lfloor X/B\rfloor$ gives a diagonal cost $O(B/X)$ in the
squared normalized sum (if $B>X$ the asserted saving costs at most
$O(\eta)$). At a surviving prime the phase difference divided by $Bh$
reduces modulo $p$ to the derivative, a monomial of degree $v-1$.
The conductor is $q_0/(h,q_0)$ with $q_0=q/B$. At dropped primes
the excluded residue classes can double, harmlessly for a fixed
number of steps. A 1-Lipschitz saving profile loses at most the
factor $(h,q_0)X^{O(\eta)}$ in a correlation. Its average loss is bounded
because $H^{-1}\sum_{h\le H}(h,q_0)\le\tau(q_0)\le(k+1)^{R_0}$.
Consequently one step acts by
\[
 g_{i+1}(u)=\tfrac12\sup_{\theta\le b\le\min(1,u)}
             \min(1-b,g_i(u-b)),
\]
with value zero for an empty interval. These profiles are 1-Lipschitz:
when $u$ decreases, clip $b$ to the new upper endpoint, and use
$g_i(0)=0$ at the boundary. To save $\sigma$ after $j$ steps one may
choose the successive shifts in
$[\theta,1-2^i\sigma]$, $1\le i\le j$, and the final
conductor in $[2^{j+1}\sigma,2-2^{j+1}\sigma]$. The sum of these
intervals is precisely the stated interval. The terminal phase has
degree at least two.
\end{proof}
\subsection{Fourier transformation at prime powers}
We give the radical dependence in the arithmetic $A$- and $B$-processes.
The $p$-adic exponent pairs of Mili\'cevi\'c \cite{Mili} treat a
prime-power modulus. Here the period of the mask can exceed the conductor.
All assertions are uniform in arbitrary additional real linear twists.

\begin{proposition}\label{prop:pairs}
Start with $(\kappa,\lambda,D)=(1/2,1/2,1/2)$ and apply a finite word
in the following transformations:
\begin{align*}
 A(\kappa,\lambda,D)
 &=\left(\frac{\kappa}{2+2\kappa},
 \frac{\kappa+\lambda+1}{2+2\kappa},
 \max\{1/2,(D+1)/2\}\right),\\
 B(\kappa,\lambda,D)&=(\lambda-1/2,\kappa+1/2,D).
\end{align*}
Read the word backwards starting at $v=k$, replacing $v$ by $v-1$ at
$A$ and by $v/(v-1)$ at $B$. Suppose that each encountered value
(including the terminal one) avoids $0,1$. Then \eqref{eq:gram-entry}
has saving
\begin{equation}\label{eq:pair-score}
 1-\lambda-\kappa(U-1)-(1+D)\theta-\eps .
\end{equation}
In particular $D\le1$ at arbitrary word length.
\end{proposition}
\begin{proof}
First restrict $n$ to one residue class modulo $m$. Locally the phase
modulo $p^e\Vert q$ is $\e(\Phi_p(n)/p^e)$, with CRT unit factors
included in $\Phi_p$. At a representative $c$ of the class, require an
integral power series on $c+p\Z_p$ whose Taylor coefficients in $t=n-c$,
modulo $p$ through order $N$, agree with
\[
 A_0(1+B_0 t)^v+C_0+D_0t,\qquad A_0B_0\in\Z_p^\times .
\]
Use $N=3+$ the length of the remaining word. For fixed $k$ the finitely
many primes at which a denominator, factorial or required coefficient is
not a unit may be excluded. Initially $v=k$ and one takes $B_0=1/c$.
Changing the representative of a residue changes none of its finite
jets modulo $p$. If $e=1$ the original class contains only one residue
modulo the conductor, so any compatible lift of that value is available.

For an interval of length at most $X$ we shall prove, when $q\ge X$,
\begin{equation}\label{eq:single}
 |T|\ll (2mqX)^\eps m^D X^{\lambda-\kappa}q^\kappa .
\end{equation}
The periodic factor has period $Q=\lcm(q,m)\le mq$. Write its unnormalised
local Fourier coefficient using the phase $\Phi(x)+b x$. For $e\ge2$
it vanishes unless $b\equiv-\Phi'(c)\pmod p$. In the remaining disk
there is a unique stationary point $x_b\in c+p\Z_p$, by Hensel's lemma
and $\Phi''(c)\in\Z_p^\times$. Splitting $x$ by translations through
multiples of $p^{\lceil e/2\rceil}$ forces $x\equiv x_b\pmod{
p^{\lfloor e/2\rfloor}}$. At even $e$ every surviving point has the
stationary phase. At odd $e$ the remaining residue sum is a quadratic
Gauss sum with unit coefficient $\Phi''(x_b)/2\bmod p$. Thus the
coefficient is $p^{e/2}\zeta\e((\Phi(x_b)+b x_b)/p^e)$, with
$|\zeta|=1$ constant on the dual disk. At $e=1$, or at a prime of
$m$ absent from $q$, the Fourier coefficient of the single class is a
linear character of $b$ times a constant of modulus one.
Fourier inversion and the geometric-sum estimate therefore give
$O((X/\sqrt Q+\sqrt Q)\log(2Q))$ for the original sum, uniformly
in its linear twist. This proves the initial pair and also gives a
saving $u/2$ at $q=X^u\le X$ with prefactor $m^{1/2}$.

For $A$, combine the input pair with completion in the profile
$g(u)=u/2$ on $0\le u\le1$ and
$g(u)=\max(0,1-u/2,1-\lambda-\kappa(u-1))$ on $u\ge1$.
It is 1-Lipschitz. When $U\ge1+\lambda$, take
$b=(\lambda+\kappa(U-1))/(1+\kappa)$ so that $U-b\ge1$ and
$1-b=1-\lambda-\kappa(U-b-1)$. We need only $b\le1$.
At every $p\mid m$ use a shift modulus $B_0$ of valuation
$\max(1,\lceil b v_p(q)/U\rceil)$. Both $B_0$ and $(q,B_0)$
lie between $X^b$ and $m X^b$. If $B_0>X$ the asserted saving is
covered by $m^{1/2}$. Otherwise take shifts $B_0h$ as before. The
normalized local difference divided by the shift reduces to the
derivative of the input jet, with one less Taylor coefficient, so its
exponent is $v-1$. A unit multiple accounts for the prime-to-$p$ part
of the conductor. The error from replacing the residual log conductor by
$U-b$ costs at most $m(h,q/(q,B_0))$. Its gcd part has bounded average.
The normalized diagonal cost is $O(mX^{b-1})$ and the correlations cost
$O(m^{D+1}X^{-g(U-b)+\eps})$. Taking the square root proves the $A$
formula. For $1\le U\le1+\lambda$ completion dominates its line:
they agree at the upper endpoint and the completion line has the
steeper negative slope.

For $B$ put $q'=\prod_{p^e\Vert q,\ e\ge2}p^e$.
All the linear characters from the other primes combine into a single
permissible real twist of the dual frequency. The remaining Fourier
coefficients have modulus $\sqrt{q'}$ and lie on one class modulo
$\rad(q')$. Their nonlinear phases have the required form. Indeed,
on a stationary disk the inverse of $b=-\Phi'(x)$ is an integral power
series in $b-b_0$, by its unit linear coefficient. Modulo $p$ its jet
agrees through the retained order with the formal inverse for
$-A_0vB_0(1+B_0t)^{v-1}-D_0$. Substitution in $\Phi(x_b)+b x_b$,
or differentiation followed by integration through that finite order,
gives a unit multiple of
$(1+B_1(b-b_0))^{v/(v-1)}$ plus an affine function, with $B_1$ a unit.
All factorials in these finite jets are invertible. Unit rescalings
and the choice of CRT coordinates preserve the assertion.

Near the frequency selected by the original real twist, an interval
of $Q/X$ frequencies has Fourier weight and total variation
$O(X/Q)$ after removal of an endpoint phase. On each dyadic interval
of distance and length $L\ge Q/X$ up to $Q$, separate the two endpoint
exponentials in the geometric progression. The remaining weight,
including the factor $1/Q$, has supremum and total variation $O(1/L)$.
For $q'\ge L$ the input pair, applied to each partial interval, costs
$\ll m^D L^{\lambda-\kappa}Q^\kappa(2mQX)^\eps$. If $q'<L$ use
completion at period $q'$. Since $Q/q'\mid m$ and $Q\ge L$, this
costs $O(m^{1/2}L^{1/2}\log(2Q))$, within the same bound.
Thus the dual contribution is bounded, up to logarithms and $X^\eps$,
by $m^D\sqrt{q'} Q^\kappa L^{\lambda-\kappa-1}$.
The largest bound occurs at $L=Q/X$. As
$\sqrt{q'}Q^{\lambda-1}\le q^{\lambda-1/2}$, this proves the
$B$ formula. Dividing $\eps$ among the finitely many operations proves
\eqref{eq:single}. Finally sum over at most $m$ unit classes.
\end{proof}

\section{Fibre estimates and the moment recursion}\label{sec:recursion}
Let $k$ and the mark intervals be fixed. We say $(M,z)$ passes the
conductor test if $1<M<k$ and for some $\sigma>2z$ the inequalities in
\cref{lem:weyl,lem:diffcells,prop:pairs}, taken together, give saving
$\sigma$ throughout $1\le U\le M$ before their arbitrarily small
losses. Shrink the maximum individual prime exponent relative to $X$ so
that a fixed positive margin above $2z$ remains. All tests are on the
entire real interval.

\begin{lemma}\label{lem:fibre}
Suppose that $(M,z)$ passes the conductor test, and let $m\le X^{M/k}$ be a product of distinct primary marks.
For $H(\alpha)=\sum_{n\le X,(n,m)=1}a_n\e(\alpha n^k)$ with
$\sum|a_n|^2\le X$, uniformly for every $\beta\in\T$,
\[
 \sum_{m^k\alpha=\beta,\ |H(\alpha)|\ge X^{1-z}/2}|H(\alpha)|^2
 \ll X^2\log(2X).
\]
\end{lemma}
\begin{proof}
The $m^k$ lifts are $(b+j)/m^k$, with $b$ representing $\beta$.
Absorb $b$ into $a_n$. The Gram matrix of the remaining phases has entries
\eqref{eq:gram-entry} with reduced denominators $q\mid m^k$.
For $d\mid m^k$, $d\le X$, retain entries whose denominators divide $d$.
This matrix is block diagonal, with each block indexed by residues mod
$d$. Finite Fourier orthogonality bounds its operator norm by
$d\max_c\#\{n\le X:(n,m)=1,\ n^k\equiv c\pmod d\}\le2k^{R_0}X$.
At each prime power a unit residue has at most $k$ roots by Hensel's
lemma. M\"obius inversion turns the condition $q\mid d$ into exact
denominator conditions. The combined part with $q\le X$ consequently
has norm $O_{k,R_0}(X)$. Every remaining entry is $O(X^{1-\sigma'})$
for a fixed $\sigma'>2z$. On $J$ lifts of magnitude at least $Y$,
duality and Cauchy's inequality, applied to unit complex coefficients
aligning their phases, give
\[
 J^2Y^2\ll X^2J+X^{2-\sigma'}J^2.
\]
For $Y\ge X^{1-z}/2$ the last term is absorbed. Dyadic summation of
$J\ll X^2/Y^2$ proves the lemma.
\end{proof}

Say that $E_u$ is available if every state satisfies
\[
 \int_{\T}|F_I(\alpha;X)|^u\dd\alpha\ll P^\eps X^{u-k+E_u}.
\]
The windows may be chosen for a finite chain of deductions and one
prescribed final accuracy. Smooth even moments are available, and H\"older
interpolation gives the straight line between any two available orders.
\begin{proposition}\label{prop:step}
Let $r\ge k+2$, $h<r$, and suppose that $E_{r-2}$ and $E_h$
are available. If $(M,z)$ passes the test and $0<L<M$, then one may take
\begin{equation}\label{eq:step}
 E_r=\max\{0,(1-L/k)E_{r-2},E_h-(r-h)z\}.
\end{equation}
\end{proposition}
\begin{proof}
Use an inner induction on the number of remaining primary marks.
For $x=\log_P X$, their exponent sum is at least $x-2\gamma$.
When $x>2\gamma/(1-M/k)$ one can choose a subset $T$ whose products
$m$ all satisfy $X^{L/k}\le m\le X^{M/k}$: accumulate their small
exponents until passing the lower endpoint, then allow for their widths.
The positive gap $M-L$ and $x\ge\gamma$ allow both errors to be made
smaller than the gap. At the other states the trivial estimate differs
from the desired one by at most $P^{O_k(\gamma/(1-M/k))}$. Finitely
many tests permit $\gamma$ to be fixed within any required error budget.

Write $F=F_I$ locally and $G_m(\alpha)=F_{I\setminus T}(\alpha;X/m)$.
Then $F(\alpha)=\sum_m G_m(m^k\alpha)$. On $|F|<X^{1-z}$ the
bound of order $h$ gives the last branch in \eqref{eq:step}. Define
$H_m$ by keeping the roots in $F$ coprime to $m$. On the remaining
set where $|H_m|\ge|F|/2$ for every $m$, weighted Jensen with the
weights $1/m$ of sum at most one gives
\[
 |F|^r\le4\sum_m m^{r-3}|G_m(m^k\alpha)|^{r-2}|H_m(\alpha)|^2.
\]
For fixed $m$ integrate with $\beta=m^k\alpha$ and apply
\cref{lem:fibre} to the conditional sum divided by $m^k$. The result is
\[
 \ll P^\eps X^2\log(2X)\sum_m m^{r-3-k}
 (X/m)^{r-2-k+E_{r-2}}
 \ll P^{2\eps} X^{r-k+(1-L/k)E_{r-2}}.
\]
On the exceptional set some $H_m$ has magnitude below $|F|/2$.
Inclusion--exclusion on the prime factors of $m$ shows that for some
nonempty divisor product $d\mid m$,
$|F|\le 2^{R_0+1}|F_{I\setminus T'}(d^k\alpha;X/d)|$.
Indeed membership $d\mid n$ singles out precisely the primes at the
corresponding indices $T'$. Sum the bound over the distinct $d$ once
for each subset of indices. By the inner induction at the same order
$r$ and Haar invariance, its cost is at most the desired bound times
a constant multiple of $\sum_{d\ge1}d^{-(r-k+E_r)}<\infty$.
This proves the induction. There are only finitely many outer steps,
so their exponent errors can be allotted in advance.
\end{proof}
\section{The linear coefficient}\label{sec:price}
For even $r$, Wooley's smooth mean value \cite[Lemma 2.1]{W95} supplies
$E_r/k\le D$, where $D\exp(D)=\exp(1-r/k)$. Fix $0<z_0<1/2$, write
$D=z_0/(1-z_0)$ and $v_0=1-D-\log D$, and begin at
$r_0=2\lceil kv_0/2\rceil$. Then $E_{r_0}/k\le D$.
This is the 1995 mean-value estimate. The 1992 paper \cite{W92}
introduced the efficient differencing used to obtain such estimates.

For a finite menu of tested triples $(L_i,M_i,z_i)$ put
$V_*(d)=\max_i\min(z_i,L_i d/2)$. Taking $h=r-2$ in
\eqref{eq:step} advances the normalized excess by
\begin{equation}\label{eq:update}
 d_{j+1}\le\max(0,d_j-2V_*(d_j)/k),\qquad d_0=D .
\end{equation}
All menu orders used here exceed $k+2$.
For $0<d_*<D$ the first even order with $d_j\le d_*$ is at most
\begin{equation}\label{eq:price}
 r_0+2+k\int_{d_*}^D\frac{\dd t}{V_*(t)} .
\end{equation}
Indeed each step ending strictly above $d_*$ decreases the integral
from $d_*$ to the current excess by at least $2/k$, by monotonicity
of $V_*$. Termination also follows from $V_*(d_*)>0$.

Suppose $z_i$ decrease, $L_i$ increase, $L_0>2(1-z_0)$ and
$2z_n/L_n<d_*$. From $D$ to $2z_0/L_0$ use the height $z_0$. Then
use $L_i t/2$ to $2z_{i+1}/L_i$ and the height $z_{i+1}$ down
to $2z_{i+1}/L_{i+1}$. Integration yields
\begin{equation}\label{eq:menu}
 \int_{d_*}^D\frac{\dd t}{V_*(t)}
 \le \frac{D}{z_0}-\frac2{L_n}
 +\sum_{i<n}\frac2{L_i}\log\frac{z_i}{z_{i+1}} .
\end{equation}
For a finite band $K\le k\le J$ take
\[
 c=\min(7/20,21/(10\lfloor\sqrt K\rfloor)),\quad
 d_*=\frac{c}{J(\log J+8.1)},\quad a=2\lceil ck/2\rceil .
\]
At the stopping order $r$, the excess satisfies
$E_r<a\rho$, where $\rho^{-1}=k(\log k+8.02113)$. Hence $S=r+a$
is even. Thus $S/k$ is at most
\begin{equation}\label{eq:finite-price}
 2+\log\frac{1-z_0}{z_0}
  +\sum_{i<n}\frac2{L_i}\log\frac{z_i}{z_{i+1}}
  -\frac2{L_n}+c+\frac6K .
\end{equation}
The three roundings come from $r_0$, the last step, and $a$.
One may increase $r$ to $4k+4$ if necessary using $|F|\le P$.
In a finite band $(4k+4+a+1)/k\le4+7/20+7/21<6.873$.
In the infinite band below the same adjustment is smaller.

\begin{lemma}\label{lem:coeff}
For each $k\ge21$ one may choose $r,a$ as above, with
$E_r<a\rho$, $r\ge4k+4$ and
\[
 \frac{S+1}{k}\le
 \frac{360936198047037988953301793}{52515099422500000000000000}
 <6.873 .
\]
The limit superior of $(S+1)/k$ is at most the fraction in
\cref{thm:main}.
\end{lemma}
\begin{proof}
We detail the finite certificate and the treatment of its infinite tail.
The witness files, exact verifiers, and their returned bounds are in
\cite{cert}. Write $b=99999/100000$ and use $L_i=bM_i$. A row with
saving requested at $\sigma>2z_i$ is accepted if the $A$-intervals
of \cref{lem:diffcells} and degree intervals from \cref{lem:weyl}
cover $[1,M_i]$, or the line \eqref{eq:pair-score} covers it.
The intervals in the first test are sorted by lower endpoint. Starting
at 1, advance through all overlapping intervals, rounding the lower
endpoints upwards and the upper endpoints downwards. A word is checked
both for its exponent pair and for the forbidden values of its model
exponent at the degrees concerned. The finite file distributes $k$ into
successive integer bands, beginning at $K=21$, $J=\lfloor8K/5\rfloor$,
then $K'=J+1$. At each band the stored $M_i$ use scale $10^4$.
The starting height $z_0$ is a displayed rational function of the finite
pair list in the verifier (and is rounded down to $10^{-5}$). Subsequent
heights are rounded down on scale $10^{32}$ after multiplication by
$400/401$. Each row contributes its term to \eqref{eq:finite-price}.
The verifier checks the required initial and terminal inequalities. In
the elementary-only version the uniform coefficient for $S/k$ is
$6.918356527842169$. With Fourier transformations, the finite bands
through $k=180000^2$ give the fraction above for $(S+1)/k$.
For example, the upper bounds for $S/k$ in the first three bands are
$6.520271,6.568709,6.633528$. These are upward roundings.

Some finite bands reach conductors above the last individually tested
degree interval. The following continuous test, also used in the
infinite tail, supplies the rest. On $F_0\le U\le M$ take
$j=\lceil hU\rceil$, $h>1$, and write $b_1=h+1/F_0$. Then
$j(j-1)\le h b_1 M^2$. The three entries in \cref{lem:weyl}
are respectively at least
\[
 (h+2/F_0)^{-1}-b_1 M\theta,\quad
 1-1/h-\theta,\quad (h-1)F_0-b_1 M\theta .
\]
It suffices that their minimum exceed $\sigma h b_1 M^2$, with
$j\le k$. Bounds for the finite bands take $F_0=120$ and test the
explicit rational slopes in \cite{cert}. The small-conductor portion
is required in the same test at the actual radical bound $\theta$.

For $k\ge K=180000^2$ start at $z_0=107/1000$. Set
$F_0=300$, $H=\sqrt k/8$, $\tau=1/14$ and take the common radical
bound $1/(8\sqrt K)$ throughout the menu. Up to conductor $F_0$,
generate all the pairs obtained by at most fourteen $A,B$ operations.
For any fixed word the numerator and denominator of the model exponent
are integral linear functions of $k$ of determinant $\pm1$.
Their coefficient sizes, including their differences, are at most
$2^{15}$, so none of the forbidden values occurs for $k\ge K$.
Use $1+D\le2$. Starting at $z_0$ decrease heights by $1000/1001$,
again rounding down at scale $10^{32}$. For each height take the largest
certified connected interval from the pair lines and degree reductions,
on conductor scale $10^6$. A strict factor $2.00001$ in the Gram
threshold leaves room for the exponent errors. This also certifies a
positive uniform margin for all sufficiently large $k$ when the
limiting menu is first computed with $\theta=0$.

Join the last row to $z=\tau/F_0^2$ at unchanged $L=bF_0$,
then use $z=\tau/M^2$, $L=bM$ for $F_0\le M\le H$. The test above
applies at $h=3/2$ because $M\theta$ here may be bounded by
$M^2/k\le1/64$ when $\theta=M/k$. At the small conductors the
last row supplies the saving with the larger common radical bound.
For the continuous part of \eqref{eq:price} one can maximize
$\min(\tau/M^2,bMt/2)$ at $M=(2\tau/(bt))^{1/3}$ until the cap.
One needs only the finitely many choices encountered in
\eqref{eq:update} at each fixed $k$. Before clipping at $H$, the
integral from $0$ to $2\tau/(bF_0^3)$ is $6/(bF_0)$, which combines
with the final $-2/(bF_0)$ in \eqref{eq:menu}. Choose
$a=2\lceil8\sqrt k\rceil$ and
$d_*=16/(k\sqrt k(\log k+8.1))$. Clipping adds at most
\[
 \frac{16}{b\sqrt k}\log\!\left(\frac{64\tau}{b}
                  (\log k+8.1)\right)
\]
since $d_H/d_*=(64\tau/b)(\log k+8.1)$ at
$d_H=2\tau/(bH^3)$. This correction and $16/\sqrt k+7/k$ decrease for
$k\ge K$. This proves coverage of the whole infinite interval.

All comparisons use integers. In evaluating the prices, for $1\le x<5/4$
one can majorize $\log x$ by $u-u^2/2+u^3/3$, $u=x-1$. Otherwise remove
powers of $2$ and use
$2\sum_{j=0}^{15}y^{2j+1}/(2j+1)+2y^{33}/(33(1-y^2))$,
$y=(x-1)/(x+1)$, as an upper bound. The cumulative prices are rounded
upwards at scale $10^{16}$. For the limiting statement the fixed menus
have vanishing radical losses. The enumeration uses 987 nonzero-slope
pairs and 11776 height levels. Its rational ceiling is the one stated
in \cref{thm:main}. Strict Gram slack, and arbitrarily small decreases
of the finitely many shortening parameters when taking the limit,
justify passing from the limiting rows to all sufficiently large $k$.
\end{proof}
\section{Representation and exceptional sets}\label{sec:represent}
Take $k\ge21$ and $a,r,S$ from \cref{lem:coeff}. Fix a positive
error smaller than the gap $a\rho-E_r$. Write $f=f_{P^\eta}(\alpha;P)$
for a smooth sum containing all the roots in $F$.
The parameter $\eta$ is first chosen sufficiently small for the
following two results, then the mark exponents are made smaller.
Br\"udern and Wooley \cite[Theorem 1.1]{BWminor} give
$\sup_{\T\setminus\M(P)}|f|\ll P^{1-\rho}$.
Their major-arc theorem \cite[Theorem 1.3]{BWmajor} gives, for some
$\nu>0$ and $1\le Q\le P$,
\begin{equation}\label{eq:pruning}
 \int_{\M(P)\setminus\M(Q)}|f|^S\dd\alpha
 \ll P^{S-k}Q^{-\nu}.
\end{equation}
Indeed one may take $t_0=\lfloor k/2\rfloor$ and
$\omega=\omega'=1$ in that theorem since $t_0>6$ and $S>2t_0+4$.
Adding $\M(1)$ bounds the full broad major-arc integral by $O(P^{S-k})$.

Put $I=\int|F|^S$. Orthogonality and the coefficientwise inclusion
$0\le 1_{\mathrm{supp}(F)}\le1_{\A(P,P^\eta)}$ give
\begin{equation}\label{eq:orth}
 I\le\int_{\T}|f|^a|F|^r\dd\alpha .
\end{equation}
Indeed expand $S/2$ positive and $S/2$ negative factors in $I$ and
replace $a/2$ factors on each side by $f$. The mixed count increases.
There is no factor for a multiplicity of marks. On the minor arcs
the right side is $O(P^{S-k-\epsilon_0})$ for some $\epsilon_0>0$.
On the major arcs H\"older bounds it by
$O(P^{(S-k)a/S} I^{1-a/S})$. We conclude that
\begin{equation}\label{eq:sharp}
 I\ll_k P^{S-k},\qquad
 \int_{\T\setminus\M(Q)} |f|^a|F|^r\dd\alpha
 \ll P^{S-k}(P^{-\epsilon_0}+Q^{-\nu a/S}).
\end{equation}
The second inequality follows from \eqref{eq:pruning} and the first.

For $Q_0$ fixed and $P\to\infty$ one has on $\mathcal K(Q_0)$
\[
 f(a'/q+\beta)=\rho_D(1/\eta)q^{-1}S_k(q,a')V(\beta)+o(P),
\]
where $\rho_D$ is Dickman's function. This is the fixed-denominator
expansion of \cite[\S6, preceding (6.15)]{BW23}. The constant is
positive. The same expansion for $F$ with coefficient $c_P$ follows from
\cref{lem:core}. Thus the Fourier coefficient of $f^a F^r$ over this
core has main term given by the truncations of
\begin{align}
 \Si_u(n)&=\sum_{q\ge1}\sum_{(b,q)=1}
       (S_k(q,b)/q)^u\e(-bn/q),\label{eq:ss}\\
 \Ji_u(y)&=\int_\R v(t)^u\e(-yt)\dd t
       =\frac{\Gamma(1+1/k)^u}{\Gamma(u/k)}y^{u/k-1}
        \quad (0<y\le1).\label{eq:si}
\end{align}
In \eqref{eq:ss} the inner sum is over reduced residues. The real
identity follows by convolution of $k^{-1}x^{1/k-1}1_{(0,1)}(x)$.
The standard complete-sum bound
$|S_k(q,b)|/q\ll_{k,\eps}q^{-1/k+\eps}$ (see \cite[Ch.~4]{Vbook})
gives uniform absolute convergence for $u>2k$. For integers with a
primitive solution in $\Z_p$ at each $p$, the Euler product
$\Si_u(n)$ is bounded below uniformly in $n$. Indeed the large-prime
factors are uniformly close to one with summable errors, since the
absolute sum over prime powers converges uniformly. At finitely many
remaining primes one unit variable has derivative of valuation $v_p(k)$.
The classes modulo $p^{2v_p(k)+1}$ with such a solution lift with
positive density bounded below, by Hensel's lemma. This also identifies
the Euler factors by character orthogonality. The classical local
Waring theorem gives primitive solutions everywhere for $u\ge4k$
($k\ge3$), and $u\ge2k+3$ suffices if $k$ is not a power of two
\cite[\S2]{VWsurvey}. Take $u=S$, $P^k/2\le n\le P^k$, choose
$Q_0$ sufficiently large in \eqref{eq:sharp}, and then let $P$ grow.
The Fourier coefficient is positive. This proves the large-degree
part of \cref{thm:main}.

For the density assertion put $t=S/2$ and use the Fourier coefficients
of $f^{a/2}F^{r/2}$. By Bessel's inequality and \eqref{eq:sharp},
the sum of squares of its Fourier coefficients outside
$\mathcal K(P^\delta)$ is $O(P^{2t-k-\epsilon_1})$ for some
$\epsilon_1>0$. The integral inside this growing core is
$\gg P^{t-k}$ for each primitive locally soluble $n$ in the indicated
dyadic interval. Here are the details needed to keep a power saving
in the exceptional set. On a fixed rectangular core both sums have
the expansions above. Between that core and $\mathcal K(P^\delta)$
bound $f$ trivially. Use \cref{lem:core} on $F$. After decreasing
$\delta$, the terms involving the approximation error contribute
$O(P^{t-k-\epsilon_2})$: the total arc length is
$O(P^{-k+3\delta})$, and the relative error in one factor is
$O(P^{-\xi/2+2\delta})$. The absolute tails of its main approximation
are bounded by the tails of
\[
 \sum_{q\ge1}q^{1-r/(2k)+\eps}
 \int_{\R}(1+|b|)^{-r/(2k)}\dd b ,
\]
which converge since $r/2\ge2k+2$. First take the fixed core large
relative to the positive local lower bound. Bessel now gives
$O(P^{k-\epsilon_1})$ exceptions on the dyadic interval. Increasing
$t$ by adding marked factors preserves the result, since their extra
powers on the minor arcs are bounded by $P$ apiece and the real
and local factors are evaluated at the new order. Summing dyadic
intervals gives the asserted power saving.

For $k=2^j\ge8$, odd $k$th powers fill $1+4k\Z_2$.
For example this follows by applying the $2$-adic logarithm on
$1+4\Z_2$, and using $\Z_2^\times=\{\pm1\}(1+4\Z_2)$.
Even $k$th powers are multiples of $2^k$. Hence for $t<4k$ a primitive
dyadic solution exists exactly on the residues $1,\ldots,t\bmod4k$.
For $t\ge4k$ all residues occur. At odd $p$ one may fix a unit variable
to be $1$. Cauchy--Davenport applied to the remaining powers, whose
residues including zero have size $1+(p-1)/(k,p-1)$, gives all residues
as soon as $t-1\ge k$. Hensel's lemma then applies. If every root is
even one may divide $n$ by $2^k$ and start again. Thus the represented
integers consist, up to the power-saving exceptions, of the disjoint
sets $2^{kj}m$ with $m\bmod4k\in\{1,\ldots,t\}$. Summing their
geometric densities and their errors proves \eqref{eq:dyad-density}.
The progression $4k\bmod8k$ excludes each smaller $t$ from density
one: its members require at least $4k$ odd summands because
$8k\mid2^k$. In degree four the familiar value is $G_1(4)=15$
\cite{VWsurvey}.

\section{Intermediate powers}\label{sec:small}
For fourth powers there are two useful routes back to Davenport's
$G(4)=16$ \cite{DavenportFour}. One follows the diminishing ranges of
Davenport and Erd\H os, combining two blocks of four biquadrates with
eight unrestricted ones. The other uses the quartic Vinogradov mean
value theorem and Hua's eighth moment. Translation and sampling control
the minor arcs, and a smooth cutoff permits Poisson summation on the
major arcs.

For seventh powers a small change in Wooley's argument \cite{W16}
recovers his bound of $31$. Take six ordinary roots and twenty-five
smooth ones. H\"older's inequality uses shares $3/28$ and $25/28$
of the ordinary minor-arc moment of order $56$ and the smooth moment
of order $28$. The resulting minor-arc exponent is
$24-60703/11200000+\varepsilon$, strictly below the main-term exponent
when $\varepsilon$ is small enough.

The November working manuscript used eight ordinary roots and thirty
smooth roots for $G(8)\le38$, and eleven ordinary and thirty-five
smooth roots for $G(9)\le46$. Both deductions use refined sampling of
the coefficients in efficient differencing. For eighth powers one also
averages a quadratic second moment according to rational height.
The corresponding rows of \cref{tab:input} below instead use seven
ordinary roots, together with thirty-one or thirty-seven marked roots.

We give the inputs and both final inequalities for every numerical
deduction at $8\le k\le20$. Let
$f_0(\alpha)=\sum_{x\le P}\e(\alpha x^k)$ and $w=k(k+1)$.
On the ordinary minor arcs one has
\begin{equation}\label{eq:ordinary}
 \int_{\m}|f_0|^w\dd\alpha\ll P^{w-k-1+\eps}.
\end{equation}
For $k\le16$ use \cite[Lemma 3.1]{W16}. In all degrees the estimate
follows by inserting the mean-value theorem at $w/2$ in
\cite[Theorem 2.1, with corrigendum]{W12}. In detail the minor-arc
bound multiplies $J_{w/2,k}(P)$ by
$P^{k(k-1)/2-1+\eps}$. The main conjecture gives
$J_{w/2,k}(P)\ll P^{w/2+\eps}$. One may take the broad ordinary arcs
with $q\le P$ and $|q\alpha-b|\le P^{1-k}/(2k)$ as in
\cite{V89,W16}; their complement lies in the minor arcs required
in \eqref{eq:ordinary}. For pruning we use
$\int_{\M_{\rm ord}(P)\setminus\mathcal K(Q)}|f_0|^v
\ll P^{v-k} Q^{-\delta_v}$, $2\le Q\le P$, for the exponents $v>k+2$
below, by \cite[Lemma 5.1]{V89} and \cite[(4.1)]{W16}. A rectangular
core contains the smaller pruning arcs. Interpolation supplies real orders.

For the initial even moments, Hua's fourth-moment bound gives
$E_4=k-2$. For orders $u=6,8,\ldots$ use
$E_u=\lambda_u+k-u$ from the degree tables of \cite[\S\S10--22]{VW00}.
The even smoothly restricted moments of order $C$ have the log-free
bound $O(P^{C-k})$ at the orders in \cref{tab:input}. For $k\le16$
this is \cite[Table 2 and (4.2)]{W16}. For $17\le k\le20$ use the
respective even orders $u=96,104,112,120$. Their published excesses
satisfy
\[
 C-u>(w-C)E_u .
\]
Indeed let $U_C=\int|f_R|^C$. Coefficientwise domination replaces
$C-u$ of its factors by $f_0$. On the minor arcs use H\"older shares
$(C-u)/w,(w-C)/w,u/w$ on the integrals of
$|f_0|^w,|f_R|^u,|f_R|^C$. By \eqref{eq:ordinary} the result is
$O(P^{(C-k)(1-u/w)-\epsilon_3}U_C^{u/w})$. On the major arcs
the ordinary $C$th moment is $O(P^{C-k})$, and H\"older again absorbs
the power of $U_C$ strictly below one. Thus $U_C\ll P^{C-k}$
and by orthogonality the same is true of $F$ at order $C$.

\begin{table}[ht]
\centering\small
\begin{tabular}{rrrrrrr}
\toprule $k$&$a$&$b$&$C$&$aC/(C-b)$&lower bound for $10^6\sigma_k$&updates\\
\midrule
8&7&31&42&$294/11$&21015&10\\
9&7&37&50&$350/13$&4361&13\\
10&7&44&60&$105/4$&10296&17\\
11&7&50&68&$238/9$&153&19\\
12&7&57&76&28&4646&24\\
13&7&63&84&28&944&103\\
14&7&70&92&$322/11$&1213&31\\
15&6&77&100&$600/23$&105&4221\\
16&7&83&110&$770/27$&870&149\\
17&7&90&118&$59/2$&3003&963\\
18&8&96&128&32&3856&1064\\
19&7&103&138&$138/5$&676&1127\\
20&7&110&148&$518/19$&1410&1245\\
\bottomrule
\end{tabular}
\caption{Ordinary and marked counts and certified minor-arc savings.}
\label{tab:input}
\end{table}
Interpolate the seed excesses by their lower convex hull. A row of
the certificate specifies $(r,h,M,z,E_r)$, the difference $M-L$, and
a strict conductor margin. Check \cref{lem:diffcells} at that margin
and \eqref{eq:step} using the previously certified orders. For instance
the first degree-eight row uses $r=25$, $h=24$, $M=199/100$,
$L=M-1/1000$, $z=5041/60000$, giving $E_r=286051517/500000000$.
All rows and the exact seed values appear in \cite{cert}. At final
order $r=bw/(w-a)$ the test is
\[
 \sigma_k=\frac aw-\left(1-\frac aw\right)E_r>0.
\]
Applying H\"older to $|f_0|^a|F|^b$ on the minor arcs gives
$O(P^{a+b-k-\sigma_k+\eps})$. On the broad major annulus apply
H\"older with shares $b/C,1-b/C$ and use the ordinary moment at
$aC/(C-b)>k+2$. Thus outside a polynomial core the integral has a
power saving. On it both $f_0$ and $F$ have the complete-sum
approximations of \cref{lem:core}, the former with coefficient one.
Since $s=a+b\ge4k$ the positive main term proves \cref{tab:small}.
The known bounds $G(1),\ldots,G(7)\le1,4,7,16,17,24,31$ complete
the proof of \cref{thm:main} \cite{VWsurvey,W16}.
For almost-all representation take
$t=\lceil a/2\rceil+\lceil b/2\rceil$. After squaring the kernel the
additional factors cost at most their sizes $P$. Bessel gives a
power saving as above, while the polynomial core has the full local
main term because $t>2k$. This completes \cref{thm:density}.
\section{Bilinear sums and almost-prime roots}\label{sec:bilinear}
The free cofactor in \cref{sec:roots} permits an arbitrary bounded
weight in each of the two factors. This provides a pointwise estimate
down to logarithmically narrow arcs, useful when the root set has
density of order $1/\log P$.
\begin{lemma}\label{lem:type2}
Let $D,Y\le P$, $|a_d|,|b_y|\le1$, and on intervals in
$[D,2D],[Y,2Y]$ put $T(\alpha)=\sum_{d,y}a_db_y\e(\alpha d^ky^k)$.
If $|\alpha-b/q|\le q^{-2}$, $(b,q)=1$, and $q\le P^k$, then
for some $C_k,c_k>0$
\[
 |T(\alpha)|\ll_k DY(\log(2P))^{C_k}
 \bigl(D^{-1}+Y^{-1}+q^{-1}+q/(DY)^k\bigr)^{c_k}.
\]
The same estimate, with an additional logarithmic factor and
$O(P^{-2})$, holds with the constraint $dy\le P$.
\end{lemma}
\begin{proof}
Cauchy's inequality in $d$ removes $a_d$. On expanding the square,
a second Cauchy inequality in $(y_1,y_2)$ removes the products
$b_{y_1}\overline{b_{y_2}}$. The fourth power of the normalized
sum is bounded by a constant times the box average of
\[
 \e\bigl(\alpha(d_1^k-d_2^k)(y_1^k-y_2^k)\bigr).
\]
This average is real and nonnegative before reparametrization. Set
$h=d_1-d_2$, $j=y_1-y_2$ and take absolute values of the interior sums.
Apply $k-1$ further Weyl differences in $d$ and $k-2$ in $y$. At each
step square the preceding average and use Cauchy's inequality in the
already fixed shifts, followed by the identity obtained by opening the
square in the chosen variable. The intersection of every translated
interval is an interval. The original shift $h$ and its $k-1$ successors
have size $O(D)$, and $j$ and its $k-2$ successors have size $O(Y)$.
A zero shift has relative cost $O_k(D^{-1}+Y^{-1})$. After the
differences the $d$-polynomial is constant and the $y$-polynomial is
linear with coefficient $(k!)^2\alpha$ times the product of these
$2k-1$ shifts. This follows by taking the difference of the leading
term at each step. Write $M=D^kY^{k-1}$. Dividing the shifting sums
by their lengths and summing the remaining progression in $y$ gives
\[
 (|T|/(DY))^{2^{2k-1}}\ll_k D^{-1}+Y^{-1}+
 \frac1M\sum_{n\ll_k M}\tau_{2k-1}(n)
       \min\!\left(1,\frac1{Y\|(k!)^2\alpha n\|}\right).
\]
The divisor function counts all ordered absolute shifts with the
indicated product. The sum of $\tau_j(n)^2$ up to $x$ is
$O_j(x(\log(2x))^{j^2-1})$. For instance
$\tau_j(p^e)^2\le\tau_{j^2}(p^e)$ by realizing two prescribed lists of
margins as a $j\times j$ matrix of nonnegative entries of sum $e$,
and the elementary divisor mean suffices. Cauchy's inequality removes
the divisor weight at the price of a square root and a fixed log
power. For the remaining mean use
\[
 \sum_{n\le Z}\min(Y,\|\alpha n\|^{-1})
 \ll (Z/q+1)(Y+q\log(2q)).
\]
This follows by separation of consecutive blocks of length at most
$q/2$ and also applies to the subsequence $(k!)^2n$ by enlargement of
$Z$. As $M\ge D$ and $MY=(DY)^k$, this proves the first assertion.
For the cutoff use truncated Perron on
$(\lfloor P\rfloor+1/2)/(dy)$ at real part $1/\log(2P)$ and height
$P^8$. Its separated factors have bounded modulus after including
the numerator and the integral costs $O(\log P)$. The distance of
$dy$ from the half-integer is at least $1/2$, giving total truncation
error $O(P^{-2})$ even before restriction to dyadic pairs.
\end{proof}

For every prescribed $A>0$ there is consequently $B=B(A,k)>0$
such that the marked sum with an arbitrary weight bounded by one on
$d$ and one on $y$ satisfies
\begin{equation}\label{eq:logminor}
 \sup_{\alpha\notin\M((\log P)^B)}|F_{a,b}(\alpha)|
 \ll_A P(\log P)^{-A}.
\end{equation}
To check the ranges, partition dyadically. Boxes with
$DY<P(\log P)^{-A-4}$ have negligible total, using the trivial
bound and the number $O((\log P)^2)$ of boxes. In the remaining boxes
both $D$ and $Y$ exceed a fixed power of $P$ by the support of $d$
and by $d\le P^{1-\xi/2}$. Dirichlet approximation supplies
$(\log P)^B<q\le P^k/(\log P)^B$ with
$|q\alpha-b|\le (\log P)^B/P^k$ (integer rounding changes only an
absolute factor). Apply \cref{lem:type2}, taking $B$ sufficiently
large compared with $A,k,C_k,1/c_k$. The same reasoning with power
cutoffs gives
\begin{equation}\label{eq:polyminor}
 \sup_{\alpha\notin\M(P^\delta)}|F_{a,b}(\alpha)|\ll P^{1-\epsilon_4}
\end{equation}
for each sufficiently small fixed $\delta>0$, by discarding boxes
with $DY<P^{1-\delta/(2k)}$.

\begin{proof}[Proof of \cref{thm:primes}]
Let $F_{\mathrm p}$ be the marked sum with $y$ restricted to primes
$y\ge P^{\xi/4}$. All marked factors are prime, lie in disjoint
intervals and exceed $y$. The roots are therefore squarefree and
$\ell=R_0+J+1$. The window exponents and $2\xi$ can be taken below
any prescribed positive bound. Every prime exceeds a positive power
of $P$, giving the stated bounds relative to $n\asymp P^k$ after
adjusting the constants.

Orthogonality gives $\int|F_{\mathrm p}|^S\le\int|F|^S
\ll P^{S-k}$. At $s=S+1$, \eqref{eq:logminor} bounds the integral
off $\M((\log P)^B)$ by $O(P^{s-k}(\log P)^{-A})$ with $A$ arbitrarily
large. We may enlarge these arcs to $\mathcal K((\log P)^B)$.
On this core all marked products are units modulo $q$. Define
\begin{align*}
 U(q,b)&=\phi(q)^{-1}\!\!\sum_{x\bmod q}^{\!*}\e(bx^k/q), \qquad C_P=\sum_d\frac1{d\log(P/d)}\asymp_k(\log P)^{-1},\\
 W(\beta)&=\sum_d\int_{P^{\xi/4}}^{P/d}
       \frac{\e(\beta(dy)^k)}{\log y}\dd y.
\end{align*}
The Siegel--Walfisz theorem \cite[Corollary 11.19]{MV} applies to
each progression in the free cofactor uniformly on this core:
$P^{\xi/4}\le y$, $q\le(\log P)^B$. Sum over the residue classes,
use partial summation in $y$, and use \eqref{eq:mass}. For any
$A'>0$ one obtains
\[
 F_{\mathrm p}(b/q+\beta)=U(q,b)W(\beta)
                 +O_{A',B}(P(\log P)^{-A'}).
\]
The log powers from the residue sum, the variation
$1+P^k|\beta|$ and the total arc length
$O(P^{-k}(\log P)^{3B})$ can all be covered by the choice of $A'$.
Multiplication by $d$ permutes the unit residues, so no distribution
of $d$ in arithmetic progressions is required.

Uniformly in $P$ one has
$|W(vP^{-k})|\ll P C_P(1+|v|)^{-1/k}$. Indeed the primitive
oscillatory integral on any initial subinterval in $y$ is
$O((P/d)(1+|v|)^{-1/k})$, and partial summation against the decreasing
function $1/\log y$ costs $O_\xi(1/\log P)$. On compact $v$-intervals
one has $W(vP^{-k})/(P C_P)=v_0(v)+o(1)$ with
$v_0(v)=\int_0^1\e(vx^k)\dd x$. To see this, put $x=dy/P$ and
first take $x\ge\epsilon$ fixed. There the logarithm equals
$\log(P/d)+O_\epsilon(1)$, uniformly in $d$. The remaining interval
costs $O_\xi(\epsilon)$. Consequently the normalized real factor
converges uniformly for $n/P^k\in[1/2,1]$ to \eqref{eq:si}.

The normalized unit singular series uses $U(q,b)$ in \eqref{eq:ss}.
For $p\nmid k$ its coefficients vanish when $p^e\Vert q$ with
$e\ge2$, by translation of the unit variable by $p^{e-1}$.
At $e=1$ the unit sum contributes $O_k(p^{-1/2})$ for nonzero
frequency by Weil. For $p\mid k$, $v=v_p(k)$, translation by
multiples of $p^{e-v-1}$ gives the same vanishing for $e\ge2v+2$.
Multiplicativity thus gives an absolutely convergent series at
every order used here. Its local factor is
\[
 \lim_{b\to\infty}\frac{p^b}{\phi(p^b)^s}
 \#\{(x_1,\ldots,x_s)\in((\Z/p^b\Z)^\times)^s:
                  \textstyle\sum x_i^k\equiv n\pmod{p^b}\}.
\]
For large $p$ it is $1+O_k(p^{1-s/2})$. At the finitely many
remaining primes an existing unit solution lifts with uniformly positive
density by fixing the root classes modulo $p^{2v_p(k)+1}$.
The unit-solubility condition is consequently a finite congruence
condition. It includes, for example, the class $n\equiv s$ modulo
a suitable fixed integer: at the exceptional primes use all ones,
while at large $p$ the positive local factor supplies a solution.
Hence for the stated classes the major arcs are
$\gg P^{s-k}C_P^s$. This dominates the complementary integral on
taking $A>s$. For $t=S/2+1$, Bessel and \eqref{eq:logminor} bound
the minor-arc square integral by $O_A(P^{2t-k}(\log P)^{-A})$ for any
$A$, while the same core calculation gives $\gg P^{t-k}C_P^t$ at
locally admissible targets. The root sets may be chosen independently at the different positions.
Keep their numbers of marks equal. Perturb the exponents within each
fixed-sum hyperplane so that all primary and auxiliary intervals for
different positions are disjoint, taking their widths after the
perturbation. Choose the finitely many cofactor exponents $\xi_i$ so
that the intervals $[\xi_i/4,2\xi_i]$ are disjoint and lie below
all the mark exponents. These choices preserve every strict estimate
above. The prime divisors of different roots are then disjoint. In the
integral argument apply H\"older to the $S$ remaining factors after
the pointwise estimate, using the same individual sharp moment at each
position. For Bessel the same argument applies to the doubled factors.
All normalized real factors have the same limit. This proves the
pairwise coprimality assertion and its almost-all version.
\end{proof}

\paragraph{Removal of the unit restrictions.}
We prove the last assertion of \cref{thm:primes}. Put
$\mathcal P_0=\{p:p\le4k+5\}$ and write $J_0=\#\mathcal P_0$.
For each $p\in\mathcal P_0$ choose an ordinary primitive local
solution $(z_{1,p},\ldots,z_{s,p})$ at $n$, which exists since $s\ge4k$.
Assign $p$ to the multiplier $d_i$ exactly when $z_{i,p}$ is a nonunit.
Complete each $d_i$ to a squarefree product of $J_0$ primes by using
a fixed private list of $J_0$ primes greater than $4k+5$ for variable $i$.
The private lists are pairwise disjoint. Only finitely many multiplier
vectors can occur.

The equation $\sum_i d_i^k u_i^k=n$ now has unit solutions at every
prime, with at least one $d_i$ a unit. At a small $p$, put
$v=v_p(k)$. As $k\ge2v+1$, the terms with $p\mid d_i$ automatically
match $z_{i,p}^k$ modulo $p^{2v+1}$. In the other terms start at
$u_i=z_{i,p}/d_i$. The derivative in one such variable has valuation
$v$, so Hensel's lemma supplies a unit solution. At a larger $p$, at
most one $d_i$ is divisible by $p$ and $p\nmid k$. The remaining
unit powers have residue set of size at least $(p-1)/k$. Their
$s-1$-fold sum is the whole residue field by Cauchy--Davenport.
If every multiplier is a unit first fix one term. The solutions lift
by a unit derivative.

Use the product $\prod_i F_{\mathrm p}(d_i^k\alpha)$ in the preceding
argument. Outside logarithmic arcs of cutoff $Q$, each $d_i^k\alpha$
is outside those of cutoff $Q^{1/2}$ for large $P$, by multiplying
back the scalar denominator. Apply \eqref{eq:logminor} to one factor
and H\"older to the other $S$ factors: dilation by a nonzero integer
preserves Haar measure. On the logarithmic core the $i$th complete
unit sum has numerator $b d_i^k$. The archimedean integral is
$(\prod d_i)^{-1}\Ji_s(n/P^k)$, since $n/P^k\le1$ and $d_i\ge1$.
At primes dividing the multipliers at least one unit coefficient
gives vanishing of the high-conductor terms as before. All other
prime factors are unchanged by unit dilation. The local lower bounds
are uniform over the finite collection of multipliers.
The roots $d_i x_i$ are squarefree since $x_i$ has only growing prime
factors and every $d_i$ is fixed. Thus each root has exactly
$\ell^*=\ell+J_0$ prime factors. This proves the assertion.

\section{Simultaneous diagonal equations}\label{sec:systems}
Let $h\ge1$, let $C$ be a fixed integral $h\times n$ matrix with
columns $\mathbf c_j$, and let $N_F(\mathbf b;P)$ count the solutions
of $\sum_j\mathbf c_j x_j^k=\mathbf b$ with every $x_j$ a marked
root. Rome and Yamagishi \cite[Theorems 1.4, 1.5 and Remark 1.6]{RY}
gave asymptotic formulae for unrestricted and smooth roots under a
disjoint-basis condition. Their smooth threshold is at most
$\lceil k(\log k+4.20032)\rceil+1$ bases. For the present support
one obtains the following linear threshold and a power-saving remainder.

\begin{theorem}\label{thm:matrix}
Let $k\ge21$ and $S$ be as in \cref{lem:coeff}. Suppose $C$ contains
a basis $B_0$ such that after deleting any one column of $B_0$ there
remain $S$ pairwise disjoint bases of $\mathbb Q^h$.
Then for some $\epsilon>0$, uniformly in $\mathbf b\in\Z^h$,
\[
 N_F(\mathbf b;P)=c_P^n P^{n-hk}
 \Si_C(\mathbf b)\Ji_C(\mathbf b/P^k)+O_{C,k}(P^{n-hk-\epsilon}).
\]
This includes $h(S+1)$ columns partitioned into $S+1$ bases, and
$hS+1$ columns when deleting any one column leaves $S$ bases.
\end{theorem}
\begin{proof}
Use common rectangular arcs with denominator $q\le Q=P^\delta$ and
width $QP^{-k}$ in the maximum norm. Outside these arcs at least one
linear form $\mathbf c_j\cdot\boldsymbol\alpha$, $j\in B_0$, lies
outside the corresponding scalar arcs with cutoff $Q^{1/(3h)}$.
Indeed simultaneous scalar approximations yield a common denominator
dividing the product of their denominators times $|\det B_0|$, by
Cramer's rule. The error is $O_C(Q^{1/(3h)}P^{-k})$. Reduce the common
fraction to primitive form. Use \eqref{eq:polyminor} on the exceptional
form and bound spare factors trivially. On the $S$ surviving disjoint
bases apply H\"older in blocks. Each block integral is
$\int_{\T^h}\prod_{j\in B}|F(\mathbf c_j\cdot\boldsymbol\alpha)|^S
\dd\boldsymbol\alpha=I^h$, since an integral nonsingular matrix
preserves normalized Haar measure under pushforward. Thus the
complement contributes $O(P^{n-hk-\epsilon})$.

Inside the common arcs use \cref{lem:core} on each linear form,
reducing its fraction if necessary. Finite complete-sum averaging gives
the same factor at the unreduced modulus. The total arc volume is
$O_C(Q^{2h+1}P^{-hk})$. Decrease $\delta$ to absorb all approximation
errors. Here the local factor $\Si_C$ sums, over primitive vectors
$\mathbf a\bmod q$, the products
$\prod_j q^{-1}S_k(q,\mathbf c_j\cdot\mathbf a)$ with character
$\e(-\mathbf a\cdot\mathbf b/q)$. The real factor integrates
$\prod_j v(\mathbf c_j\cdot\mathbf t)
\e(-\mathbf t\cdot\mathbf b/P^k)$ over $\R^h$.

For absolute tails, discard the factors outside $S$ disjoint bases.
Use H\"older on the primitive vectors. Within one basis write $q_i$
for the reduced denominators of its $h$ coordinates. Their least common
multiple is at least $q/|\det B|$ by the adjugate identity. The
kernel size of the map modulo $q$ is bounded in terms of $B$ (use its
Smith normal form). At given $q_i\mid q$ the number of possible images
is at most $\prod_i q_i$. Raising all $h$ factors to the $S$th power,
using the complete-sum estimate and summing over the choices of $q_i$,
costs $O_{C,\eps}(q^{1-S/k+\eps})$. The resulting tail has a power
saving in $Q$. For the real tail apply the same H\"older inequality
and change variables in each basis. Whenever $|\mathbf t|_\infty>Q$,
at least one transformed coordinate has modulus $\gg_C Q$, and
$|v(x)|^S\ll(1+|x|)^{-S/k}$. This gives a power saving again.
These estimates justify completion of the truncated factors and
prove the formula. When the real fibre has a nonsingular interior
point and there are nonsingular local points everywhere, both factors
are positive by the real and $p$-adic implicit function theorems.
\end{proof}
\section{Pointwise estimates and small fractional parts}\label{sec:fractional}
There are dense sets of smooth roots on which the saving in a pointwise
Weyl estimate has order $1/k$. The set used here has an outer group of
prime factors of product just below the square root of the full scale.
\begin{proposition}\label{prop:sampledweyl}
Fix $k\ge21$, $0<\lambda<1/2$, and $\eta>0$. Let $S$ be the order
in \cref{lem:coeff}. For every sufficiently large $P$ there is a set
$\mathcal B_P\subset\A(P,P^\eta)$ with $|\mathcal B_P|\gg P$ such
that, simultaneously for $1\le Q\le P$,
\[
 \sup_{\alpha\notin\M(Q)}
 \left|\sum_{x\in\mathcal B_P}\e(\alpha x^k)\right|
 \ll_\eps P^{1+\eps}(Q^{-1}+P^{-\lambda})^{1/S}.
\]
In particular the sum is $O(P^{1-1/(13.75k)})$ off $\M(P)$.
\end{proposition}
\begin{proof}
Choose $\lambda<\lambda'<1/2$. Take disjoint prime windows
of the form in \cref{sec:roots}, all with upper endpoint below $P^\eta$,
so every product $v$ of their marks satisfies $P^\lambda<v<P^{\lambda'}$.
Fix the positive exponents of these windows. Partition the products
into dyadic intervals $M\le v<2M$. Within an interval put $U=P/(2M)$
and let $h$ be the marked Weyl sum at scale $U$ from \eqref{eq:sharp},
with all its prime factors below the smallest outer prime. Include the
roots $vu$ for $u$ in this inner set. The outer windows recover $v$
uniquely. The number of roots is $\gg P\sum_v v^{-1}\gg P$ by
\eqref{eq:mass}. There are $O(\log P)$ intervals.

Partition the frequency circle into cells of length comparable to
$U^{-k}$. For two $v,w$ whose frequencies $\alpha v^k,\alpha w^k$
are in the same cell we have, with $d=v^k-w^k$, an integer $c$ such
that $|d|\ll M^k$ and $|\alpha d-c|\ll U^{-k}$. The determinant of any
two pairs in these ranges is $O((M/U)^k)=o(1)$, so all pairs lie on a
single rational line. If there is no nonzero pair, each cell contains
at most one sample. Otherwise the line has equation $c=ad/q$ with
$(a,q)=1$, $q\ll M^k$. Put $\beta=\alpha-a/q$.

Group the outer products in the interval by $d_0$, the product of their
primes dividing $q$. There are $O(1)$ groups since the outer primes
exceed a fixed power of $P$ and $q\ll M^k$. In one group write $v=d_0v'$.
There are $O(M/d_0)$ possibilities. Within one cell the powers $(v')^k$
agree modulo $q'=q/(q,d_0^k)$, and $v'$ is a unit modulo $q'$.
Each unit has $O_{k,\eps}((2q')^\eps)$ possible $k$th roots, by the
structure of the groups of units modulo prime powers. Moreover
$|\beta(v^k-w^k)|\ll U^{-k}$. The mean value theorem on $[M,2M]$ puts
the quotients $v'$ in an interval of length
$O(M/(d_0(1+P^k|\beta|)))$. Since $q'\ge q/d_0^k$, and
$q(1+P^k|\beta|)>Q$ off $\M(Q)$, the number of samples per cell is
\[
 B\ll P^\eps(1+M d_0^{k-1}/Q).
\]
In the case of no nonzero pair take $B=1$, $d_0=1$.

The large sieve for the polynomial $h^{S/2}$ of degree $O(U^k)$ gives
$\sum_v|h(\alpha v^k)|^S\ll B U^k\int |h|^S\ll B U^S$ in the group.
Indeed the samples split into $O(B)$ sets with separation $\gg U^{-k}$
by taking alternate cells, splitting at the cut of the circle if necessary.
H\"older and $S\ge k$ now give
\[
 \left|\sum_{v\,\text{in the group}}h(\alpha v^k)\right|^S
 \ll (M/d_0)^{S-1} B U^S
 \ll P^{S+\eps}(M^{-1}+Q^{-1}).
\]
Summing over the groups and intervals proves the estimate.
Carrying out the construction with $\max(\lambda,4999/10000)$ in place of $\lambda$ and using $S<6.873k$ leaves a positive
margin for $\eps$ in the last assertion.
\end{proof}

\begin{theorem}\label{thm:fractional}
For each $k\ge21$ and all sufficiently large $P$, uniformly for
$\alpha\in\R$,
\[
        \min_{1\le n\le P}\|\alpha n^k\|\le P^{-1/(12.31k)} .
\]
For all sufficiently large $k$ one may replace $12.31$ by $12.04$.
Here $\|x\|$ is distance to the nearest integer.
\end{theorem}
Br\"udern and Wooley \cite[Theorem 1.2]{BWminor} proved the estimate
with denominator $k(\log k+8.02113)$ for $k\ge6$. Our estimate
uses the sharp moment for a specially chosen inner set of roots.

\begin{proof}
Use the construction of \cref{prop:sampledweyl}, now with
$1/2<\lambda_-<\log_P v<\lambda_+<0.58$. On the complement
of $\M(P^{0.58})$ the calculation in its proof gives a saving
$\lambda_-/S$ on every dyadic interval for which the pairs
$(d,c)$ in the indicated rectangle have rank at most one.

Suppose that they have rank two. The lattice
\[
 \{(M^{-k}d,U^k(\alpha d-c)):d,c\in\Z\}
\]
has determinant $(U/M)^k$ and two independent vectors of bounded length.
These span a sublattice of bounded covering radius. Disjoint fundamental
Voronoi cells around points of the full lattice in a fixed box are
therefore contained in a bounded enlargement. The number of pairs
in the box is $O((M/U)^k)$, also for any fixed constants in its sides
for which the rank is two. Fix such sides large enough for the compact
Fourier support used below before distinguishing the ranks.

For integers $t,w\ge1$, let $b_y$ be the coefficients of $h^t$.
They satisfy $\sum b_y\ll U^t$ and
$\sum b_y^2\ll U^{2t-k+E_{2t}+\eps}$ by the even smooth moment
of \cite[Lemma 2.1]{W95}. The outer marks may also be taken sufficiently
small relative to $M$ for the corresponding bound at $2w$. Choose
$|a_v|\le1$ to align the phases of $h(\alpha v^k)^t$ and put
$V_y=\sum_v a_v\e(\alpha yv^k)$. Two applications of H\"older give
\[
 \left|\sum_v h(\alpha v^k)\right|^{2tw}\!
 \le (2M)^{2w(t-1)}(\textstyle\sum b_y)^{2w-2}
                (\textstyle\sum b_y^2)\sum_{y\le tU^k}|V_y|^{2w}.
\]
In the second application the three factors in the summand are
$b_y^{1-1/w}$, $b_y^{1/w}$ and $|V_y|$. Majorize the last interval
by a nonnegative Schwartz function of $y/U^k$ with compact Fourier
support, obtained as the squared modulus of the Fourier transform
of a suitably dilated smooth bump. In the expansion of $|V_y|^{2w}$
each coefficient at a power-sum difference $d$ is
$O(M^{2w-k+E_{2w}+\eps})$, by Cauchy's inequality on a translate
of the $w$-fold representation function. Poisson summation in $y$
costs $O(U^k)$ for each of the $O((M/U)^k)$ lattice pairs. Consequently
\[
 \left|\sum_v h(\alpha v^k)\right|\ll
 P^{1+\eps}\bigl(M^{E_{2w}}U^{-k+E_{2t}}\bigr)^{1/(2tw)} .
\]

Write $x=t/k$, $y=w/k$, $D_x=E_{2t}/k$ and $D_y=E_{2w}/k$.
At $M=P^\ell$, with an $O(1/\log P)$ endpoint allowance, we have
thus saved $1/k$ times the smaller of
\[
 \frac{\ell}{6.873},\qquad
 \frac{(1-\ell)(1-D_x)-\ell D_y}{2xy}.
\]
Take $t=w=k$ and $0.5588\le\ell\le0.5590$. The smooth bound gives
$D_1<0.2785$. Both entries exceed $1/12.31$. For the large-degree
conclusion take $t=\lceil3k/4\rceil$, $w=\lceil3k/2\rceil$ and
$0.57089\le\ell\le0.57093$. We have $D_x<0.40468$ and
$D_y<0.12003$, and both entries exceed $1/12.04$ as $k\to\infty$.
The three defect comparisons follow by substituting these rational
bounds in $D\exp(D)=\exp(1-2x)$; in each case the polynomial
$\sum_{j=0}^{13}(D+2x-1)^j/j!$ already exceeds $1/D$ at the
stated upper bound and the nominal $x=1,3/4,3/2$ respectively.
Choose the outer windows strictly inside the indicated ranges.
All comparisons have positive slack, which absorbs the logarithms,
the dyadic endpoint changes and the permissible moment losses.

Let $\rho$ be the claimed exponent and choose $\nu>0$ small enough
that the resulting minor bound for the dense set $\mathcal B_P$ is
$O(P^{1-\rho-\nu})$. If no $n\le P$ has
$\|\alpha n^k\|\le P^{-\rho}$, the interval $[0,P^{-\rho}]$
contains none of the fractional parts from this set.
The Erd\H{o}s--Tur\'an inequality \cite[Ch.~2, Theorem 2.5]{KN},
with $H=\lceil P^{\rho+\nu}\rceil$, then gives some $1\le h\le H$
with $|\sum_{n\in\mathcal B_P}\e(h\alpha n^k)|\gg P^{1-\rho}/\log P$.
It follows that $h\alpha\in\M(P^{0.58})$. Thus $q\le P^{0.58}$ and
$|qh\alpha-a|\le P^{0.58-k}$ for some $a\in\Z$. Choose $n=qh$.
After decreasing $\nu$, we have $n\le P$ and
\[
 \|\alpha n^k\|\le (qh)^{k-1}P^{0.58-k}<P^{-\rho},
\]
a contradiction.
\end{proof}

\section{Mixed powers and arithmetic progressions of exponents}\label{sec:mixed}
For a fixed tuple $2\le k_1\le\cdots\le k_s$ put
$H=\sum_{i=1}^s 1/k_i$. The sharp moments at the separate scales
$n^{1/k_i}$ give a bound depending on this reciprocal sum alone.
\begin{theorem}\label{thm:mixed}
If $H\ge6.873$, every sufficiently large $n$ has
$\gg n^{H-1}$ representations
\[
 n=x_1^{k_1}+\cdots+x_s^{k_s},\qquad x_i\in\mathbb N .
\]
If $H\ge3.4365$, the same lower bound holds outside
$O(N^{1-\delta})$ integers $n\le N$ with a primitive local solution
at each prime, for some $\delta>0$. When $s\ge4k_1$ the latter
statement holds with unrestricted $n$. The constants may depend on the tuple.
\end{theorem}
\begin{proof}
For each $k\ge2$ we use an even order $B_k<6.873k$ for which
the marked sum has moment $\ll P^{B_k-k+\eps}$. For $k\ge21$ use
\eqref{eq:sharp} and for $8\le k\le20$ take $B_k=2t$ in
\cref{tab:small}. Indeed replacing $2\lceil a/2\rceil$ factors in this
even moment by the ordinary sum, equally on both sides, bounds it by
the integral of $|f_0|^a|F|^b$ times $P^{2t-a-b}$. The bound follows
from \cref{sec:small}. For $2\le k\le5$ take $B_k=2^k$ by Hua's lemma
\cite[Lemma 2.5]{Vbook}. For $k=6,7$ take $B_k=38,48$.
At order $32$ Hua gives excess $k-5$ for the ordinary sum, while
\eqref{eq:ordinary} gives excess $-1$ at order $k(k+1)$ on the
minor arcs. Interpolation gives nonpositive excess at the indicated
orders. On the broad major arcs use the ordinary moment bound of
\cref{sec:small}; the marked even moments follow by orthogonality.
The Type~II estimate \eqref{eq:polyminor} applies to the marked
sums in every one of these degrees.

Take $N/2\le n\le N$, $P_i=N^{1/k_i}$, and let
$g_i=F_i/P_i$ denote the normalized marked sum of degree $k_i$.
Choose a common rectangular core $q\le N^\tau$,
$|\alpha-a/q|\le N^{-1+\tau}$ with $\tau>0$ sufficiently small.
Off this core $|g_i|\ll N^{-\delta_i}$ by \eqref{eq:polyminor},
and $\int|g_i|^{B_{k_i}}\ll N^{-1+\eps}$. For $u=1$ or $2$,
suppose $J=\sum_i u/B_{k_i}>1$. H\"older with exponents
$J B_{k_i}/u$ on the factors $|g_i|^{u/J}$ gives
\[
 \int_{\rm off\ core} \prod_i|g_i|^u\dd\alpha
 \ll N^{-1+\eps-\sum_i\delta_i u(1-1/J)} .
\]
Choose $\eps$ smaller than the fixed saving. The assumed bounds on
$H$ ensure $J>1$ at $u=1$ and $u=2$, respectively.

On the common core \cref{lem:core} replaces every $g_i$ with
$c_{P_i}q^{-1}S_{k_i}(q,a)v_{k_i}(N\beta)$, where the subscript on
$v$ denotes its degree. The normalized error after integration is
$O(N^{-1-\delta'})$ on decreasing $\tau$: the core has measure
$O(N^{-1+3\tau})$ and each factor has a power-saving error.
The series tails are bounded by $\sum_{q>Q}q^{1-H+\eps}$ and the integral
tails by $\int_Q^\infty(1+b)^{-H}\dd b$. Both permit power savings.
The resulting singular integral at $y=n/N$ equals
$y^{H-1}\prod_i\Gamma(1+1/k_i)/\Gamma(H)$, by convolution.
The mixed singular series is bounded below on the primitive classes
by the same Hensel argument as for \eqref{eq:ss}; the maximum of
$v_p(k_i)$ over this fixed tuple suffices. For $H>2$ and $s\ge4k_1$,
it is uniformly bounded below on every integer by
\cite[Lemma 7.2]{BWmix}. At $H\ge6.873$ the latter summand condition
is automatic. The off-core $L^1$ estimate proves the first assertion.
Bessel applied to the off-core $L^2$ estimate leaves $O(N^{1-\delta})$
exceptions on a dyadic interval, giving the second.
\end{proof}

Let $R(k,d)$ be the least $s$ such that the degrees
$k,k+d,\ldots,k+(s-1)d$ represent all sufficiently large integers.
Br\"udern and Wooley obtained $R(k,1)\le100(k+1)^2$ and
$R(k,k)\le(6k+6)^{2k}$, and discussed the expected linear order
for consecutive exponents \cite[Theorem 1.3 and Corollaries 1.4
and 1.5]{BWmix}.
\begin{corollary}\label{cor:mixed-prog}
For all $k\ge2$ and $d\ge1$ one has
\[
 R(k,d)\le\left\lceil\frac{k}{d}(\exp(6.873d)-1)\right\rceil .
\]
In particular $R(k,1)\asymp k$. With $C_\infty$ as in \cref{thm:main},
\[
 \limsup_{k\to\infty}\frac{R(k,1)}{k}\le e^{C_\infty}-1<964,
 \qquad
 \limsup_{k\to\infty}\frac{\log R(k,k)}{k}\le C_\infty .
\]
\end{corollary}
\begin{proof}
The reciprocal sum is at least
$d^{-1}\log(1+ds/k)$ by comparison with the integral. For tuples with
least degree tending to infinity one may replace $6.873$ by any
$C>C_\infty$, since $B_j/j<C$ for every sufficiently large $j$.
For the lower bound, a tuple with $H<1$ represents at most
$N^H$ integers up to $N$. In particular
$\liminf_{k\to\infty}R(k,1)/k\ge e-1$, by the reverse integral
comparison with error $O(1/k)$.
\end{proof}
\section{Three positive cubes and the six-variable cubic}\label{c3:density}
Write $f(\mathbf{x})=x_1^3+x_2^3+x_3^3$, and let $R(A)$ count the integers
in $[1,A]$ represented by $f$ with strictly positive integral arguments.
Maynard \cite{MaynardCubes} surveys the progress from Davenport to the smooth Weyl sum estimates of Vaughan and Wooley. Wooley obtained $R(A)\gg A^{0.91709477}$ \cite{c3:wooley}.
Three figures will be useful in what follows. The upper and lower
densities are respectively the limit superior and limit inferior of
$R(A)/A$.
\begin{center}
\small
\begin{tabular}{lll}
\toprule
Quantity & Value & Basis \\
\midrule
Upper density & $<0.11683$ & Congruence count, \eqref{c3:upper} \\
Predicted density & $0.0999425\ldots$ & Deshouillers, Hennecart and Landreau \cite{c3:dhl} \\
Lower density & $>0.09521$ & Full weighted Hooley--Manin asymptotic, \eqref{c3:surplus-bound} \\
\bottomrule
\end{tabular}
\end{center}
The elementary upper density bound is
\[
 \lambda=\Gamma(4/3)^3/6=0.1186788\ldots,
\]
the volume of the strictly ordered chamber
$D=\{0<x_1<x_2<x_3,\ f(\mathbf{x})<1\}$.
Triples with two coordinates equal contribute $O(A^{2/3})$ targets.

A congruence refinement gives the unconditional bound
\begin{equation}\label{c3:upper}
 \limsup_{A\to\infty}R(A)/A<0.11683 .
\end{equation}
Indeed, let $n_m(r)$ count triples modulo $m$ with sum of cubes
$r\bmod m$. The ordered chamber contains
$(\lambda n_m(r)/m^3)A+o_m(A)$ triples in that class, while there are
$A/m+O(1)$ available targets. Taking the smaller of the two bounds in
each class gives
\[
 \limsup R(A)/A\le \frac1m\sum_{r\bmod m}
          \min\!\left(1,\frac{\lambda n_m(r)}{m^2}\right).
\]
For reproducibility take $m=729\cdot512\cdot343\cdot13\cdot19\cdot31
\cdot37\cdot43$. For each displayed factor $b$ cube the residue-count
vector $(\#\{x\bmod b:x^3=j\})_{j\bmod b}$ by cyclic convolution.
Count the frequencies of its entries to obtain the distribution of
$n_b(r)$. By the Chinese remainder theorem multiply one entry from each
distribution, with multiplicities multiplied, to obtain the distribution
of $n_m(r)$. With $11869/100000$ in place of $\lambda$, the displayed
upper bound is the rational number
\[
 \frac{728028583078263256286286854758069}
      {6231646102240614790247825894400000}<0.11683 .
\]
The gamma estimate at the end of the proof below verifies the upward rounding.

Consider $F(\mathbf{x},\mathbf{y})=f(\mathbf{x})-f(\mathbf{y})$.
Over $\mathbb Q$ the change $\mathbf y\mapsto-\mathbf y$ identifies its
projective zero set $V$ with the Fermat cubic fourfold. Hooley and
Heath-Brown \cite{c3:hooley,c3:hooleyk,c3:hb} obtained an upper bound
$O_\varepsilon(P^{3+\varepsilon})$ for the number of integral zeros in the
box of radius $P$, conditionally on Hasse--Weil hypotheses for the smooth
hyperplane sections $V_{\mathbf c}=V\cap\{\mathbf c\cdot\mathbf X=0\}$.
Cauchy's inequality then gives $R(A)\gg_\varepsilon A^{1-\varepsilon}$.
Wang \cite[Theorem 1.7]{c3:wang} obtained $O(P^3)$, and hence positive
lower density, under three explicit conjectures. His hypothesis HW2 concerns
automorphy and GRH for the normalized middle-cohomology $L$-functions of
the smooth $V_{\mathbf c}$, their tensor, symmetric and exterior squares,
the ambient $L(s,V)$, and $\zeta(s)$. His R2 estimates a two-point
mean of reciprocal $L$-functions across the hyperplane sections.
His SF controls large prime square divisors of the dual form $F^\vee(\mathbf c)$.
In particular, for $P_0\le Z^{3/2}$ it asks for a power saving in $P_0$ in
the number of $\mathbf c\in[-Z,Z]^6\cap\mathbb Z^6$ for which
$p^2\mid F^\vee(\mathbf c)$ for some $P_0\le p<2P_0$.
For the positive-density theorem the log-free, smoothed mean-square
upper bound R2$'$ of \cite[Conjecture 6.4 and Remark 6.5]{c3:wang}
already suffices in place of the asymptotic R2. It asks for square-root
cancellation on average in smoothed sums of the reciprocal coefficients
at moduli $n$ of size up to $Z^3$. This permits an approach by moment
upper bounds for families of $L$-functions without determining the full
two-point main term. The separate density-one prediction in the classes $a\not\equiv\pm4\pmod9$
allows signed cubes.

The weighted asymptotic conjectured for $F$ has two main terms, the usual
Hardy--Littlewood main term and the contribution of the rational
three-dimensional vector spaces on $F=0$, with error $o_w(P^3)$
\cite[(1.2)--(1.3)]{c3:wang}. Restricting this prediction to a single
ordered chamber gives a numerical lower bound close to the conjectural density.

\begin{proposition}\label{c3:quant}
Assume the weighted Hooley--Manin asymptotic just stated for
$F=f(\mathbf{x})-f(\mathbf{y})$. For its singular series
$\mathfrak S=\prod_p\sigma_p$, one has
\[
 \liminf_{A\to\infty}\frac{R(A)}{A}
 \ \ge\ \lambda-\tfrac12\mathfrak S\lambda^2\ >\ 0.095.
\]
\end{proposition}

\begin{proof}
Choose $0\le w\le1$ smooth and compactly supported in $D$. A triple
$\mathbf a$ receives weight $w(\mathbf a/P)$, where $P=A^{1/3}$.
For each target the sum of its weights less the sum of the pairwise
products over unordered pairs of distinct triples is at most the indicator
of representability. This is the two-term Bonferroni inequality, applied
to independent events with probabilities given by the weights.

Only the identity plane $\mathbf x=\mathbf y$ meets $D\times D$.
Indeed, the rational planes come from pairings
\cite[Proposition 3.6]{c3:planes}. Positivity requires each pair to join
opposite sides, and strict ordering determines the permutation. Subtracting
the identical triples therefore removes the full plane term, even for a
nonconstant weight. If $g_w(t)\,dt$ is the pushforward of
$w(\mathbf x)\,d\mathbf x$ by $f$, the singular integral for the product
weight is $\int g_w(t)^2\,dt$. The asserted asymptotic and the lattice
point estimate for the first moment yield
\[
 \liminf R(A)/A\ge \int_D w-\tfrac12\mathfrak S\int_0^1g_w(t)^2\,dt.
\]
Homogeneity gives $\operatorname{vol}\{\mathbf x\in D:f(\mathbf x)<t\}
=\lambda t$ for $0\le t\le1$. Smooth cutoffs increasing to $1_D$
therefore have $g_w(t)\uparrow\lambda$ almost everywhere in $(0,1)$.
The displayed lower bound follows by monotone convergence.

For the numerical inequality write $M_p=\#\{F=0\pmod p\}$. When $p\ne3$,
a nonzero solution modulo $p$ is smooth, and a vector divisible by $p$
contributes $p^{-3}\sigma_p$ to the density. Consequently
\[
 \sigma_p=\frac{M_p-1}{p^5(1-p^{-3})}.
\]
At $3$ a primitive vector has at least one unit coordinate. On any unit
class $r+3\mathbb Z_3$ the cube map is uniform onto $r^3+9\mathbb Z_3$,
since $( (r+3u)^3-r^3)/9$ has unit derivative and is a permutation
modulo $3$. Cubes modulo $9$ take the values $0,1,-1$ with equal mass.
The constant coefficient of $(1+z+z^{-1})^6$ is $141$, including the
all-zero residue vector. Thus
$\sigma_3=9(141-1)3^{-6}/(1-3^{-3})=70/39$.

For clarity, the following finite computation uses only integers. For
each prime $p\le2000$ other than $3$, form the vector
$b_j=\#\{u\bmod p:u^3=j\}$ and its cyclic convolution
$t=b*b*b$. Then $M_p=\sum_j t_j^2$. Multiplying the displayed rational
factors gives
\[
 3.34949<\prod_{p\le2000}\sigma_p<3.34950.
\]
For $p\equiv2\pmod3$, cubing is bijective modulo $p$ and $M_p=p^5$.
For $p\equiv1\pmod3$ let $S(a)=\sum_{x\bmod p}\exp(2\pi i ax^3/p)$.
For $a\ne0$ the two cubic Gauss sums give $|S(a)|\le2\sqrt p$, and
orthogonality gives $\sum_{a\ne0}|S(a)|^2=2p(p-1)$.
It follows that $M_p=p^5+p^{-1}\sum_{a\ne0}|S(a)|^6\le p^5+32p^3$.
The primes contributing this last term are $1\pmod6$. With $M=2000$,
comparison with integrals now gives
\[
 \sum_{p>M}\log\sigma_p
 \le 32\bigl(M^{-2}+(6M)^{-1}\bigr)+M^{-2}.
\]
Hence $\mathfrak S<3.36$. With $m=10^4$ and
$U=\tfrac m6\prod_{j=1}^{m-1}(3j/(3j+1))^3$, log-convexity of the gamma
function gives $U(m/(m+1/3))^2\le\lambda\le U$, so that
$0.11867<\lambda<0.11869$.
The function $x-1.68x^2$ increases on this interval, and its value at
$0.11867$ exceeds $0.095$.
\end{proof}

This calculation needs an asymptotic only on products of smooth cutoffs
inside the positive chamber. Even a one-sided estimate on the distinct
pairs there, with the same leading constant, would suffice.

\subsection{The contribution forced by congruences}
\label{c3:congruence-surplus}
Congruences improve the conditional lower bound as well as the
unconditional upper bound. Let $g$ be the continuous piecewise linear
function on $[0,\infty)$ taking the values $0$ at $0,1,2$ and
$\binom{j-1}{2}$ at integers $j\ge3$. It is convex.
Put $D_m(a)=n_m(a)/m^2$ in the notation of
\eqref{c3:upper}, and give the residue classes modulo $m$ equal weight.
Under the hypothesis of Proposition~\ref{c3:quant} one has
\begin{equation}\label{c3:surplus-bound}
 \liminf_{A\to\infty}\frac{R(A)}{A}
 \ \ge\ \lambda-\tfrac12\mathfrak S\lambda^2+
                       \frac1m\sum_{a\bmod m}g(\lambda D_m(a))
 \ >\ 0.09521
\end{equation}
for the modulus used in \eqref{c3:upper}. No arithmetic-progression
strengthening of the conjectured asymptotic is needed.

Indeed, for a nonnegative integer $j$ one has
$\mathbf1_{j>0}=j-\binom j2+g(j)$. With a smooth cutoff $w$ as in
the proposition, retain each triple independently with its weight
as probability. Apply the identity to the number of triples retained
at a target, and then use Jensen's inequality on $g$. The probability
that a retained triple exists is at most the indicator of representability.
Within a fixed class $a\bmod m$, a second application of Jensen gives
a contribution per target at least $g(D_m(a)\int_D w)+o_{m,w}(1)$
from this last term. In fact counting the triples in the finitely
many residue boxes gives $(\int_D w)n_m(a)A/m^3+o_{m,w}(A)$
for their first moment. The pair term is still estimated by the
original conjecture over the whole chamber. Let the cutoffs increase
to its indicator as before. In the finite computation,
$g(u)=(j-1)(u-(j+2)/2)$ for $j=\lfloor u\rfloor\ge2$.
The same cyclic convolutions used for the upper estimate give
\[
 \frac1m\sum_{a\bmod m}g\bigl(0.11867D_m(a)\bigr)
 =\frac{313405679203461203128382573611121}
        {1512212787477055855766805750374400000}
 >0.00020724 .
\]
Using $\mathfrak S<3.36$ gives \eqref{c3:surplus-bound}.

There is a common limit for these congruence refinements. Give
$\widehat{\mathbb Z}$ Haar probability and take $m_j=j!$.
The nonnegative variables $D_{m_j}$ form a martingale: averaging over
the residue classes above any fixed class restores its triple count
with the appropriate volume factor. Their second moments are
\[
 \mathbb E D_m^2=m^{-5}
 \#\{(\mathbf x,\mathbf y)\bmod m:f(\mathbf x)=f(\mathbf y)\}.
\]
Along $j!$ these increase to $\mathfrak S$. For each prime the
local mean squares increase from one to the local factor by conditional
expectation, and the Chinese remainder theorem gives the products;
the convergent product and its tail bound were proved above.
Consequently $D_{m_j}\to D_\infty$ almost surely and in $L^2$, with
$\mathbb E D_\infty=1$ and $\mathbb E D_\infty^2=\mathfrak S$.
The upper estimates decrease to
$\mathbb E\min(1,\lambda D_\infty)$, and the surplus terms
increase to $\mathbb E g(\lambda D_\infty)$.
For the latter convergence use the bound
$|g(u)-g(v)|\le(1+u+v)|u-v|$ and Cauchy's inequality.
Every finite modulus divides some $j!$, so these limits give
the full refinements obtainable from the two displayed congruence
inequalities.

\subsection{Targets with many primitive representations}
Hooley \cite[Theorem 3]{c3:hooleyk} obtained a lower bound of
order $P^3$ for nontrivial equal sums of three cubes. Congruence
conditioning also gives control of the individual roots and the
multiplicities of the targets.
\begin{proposition}\label{c3:manyprimitive}
Fix $j\ge2$ and three nonempty bounded open intervals
$I_1<I_2<I_3$ in $(0,\infty)$. For all sufficiently large $P$,
at least $c P^3$ triples $x_i\in P I_i$ belong to fibers of
$f$ with at least $j$ triples in these same intervals, each
consisting of pairwise coprime squarefree integers. Here and below
these triples alone are counted, and $c>0$ may depend on the
intervals and on $j$. The chamber therefore contains
$\gg P^3$ primitive rational points off the planes on $f(\mathbf x)
=f(\mathbf y)$, with six distinct roots.
\end{proposition}
\begin{proof}
For a prime $p\equiv1\pmod3$ the normalized count of
$f(\mathbf x)=1\pmod p$ is
\[
 p^{-2}\#\{\mathbf x\bmod p:f(\mathbf x)=1\}
                    =1+6/p+O(p^{-3/2}).
\]
Indeed, with a cubic Gauss sum $g$ and a cubic character $\chi$,
$S(a)=g\bar\chi(a)+\bar g\chi(a)$ satisfies
$S(a)^3=g^3+\bar g^3+3pS(a)$ for $a\ne0$.
Orthogonality gives the formula since one has three nonzero cube
roots of one. The normalized count on imposing local squarefreeness
and pairwise coprimality is still $1+6/p+O(p^{-3/2})$.
Two zero coordinates cost $O(1)$ residue vectors at level one.
There are only $O(p)$ with a specified zero coordinate, and just
$1/p$ of the lifts of that coordinate vanish modulo $p^2$.

Take a finite set of such primes, of product $m$, and require
$f(\mathbf x)=1\pmod m$. At every other prime $q$ the density of
the two root restrictions is
$(1-1/q)^3(1+3/q)=1+O(q^{-2})$. The infinite product is positive.
The elementary sieve gives the product of all these densities in
the root box. Indeed failure at an unselected prime $q>L$ requires
either $q^2\mid x_i$ or $q\mid x_i,x_k$ for $i\ne k$.
Summation of the respective lattice bounds
$O(P^3/q^2)$ gives $o_{L\to\infty}(P^3)$ uniformly for large $P$;
primes larger than the permitted root or its square root cannot
occur in the corresponding sum. Having first imposed the bounded
conditions by the Chinese remainder theorem, let $L$ tend to infinity.

All resulting targets lie in an interval of length $O_I(P^3)$, with
at most $C_I P^3/m+O(1)$ available targets in the chosen progression.
The ratio of the count of triples to this quantity can be made
arbitrarily large: the product of the gains at selected primes
diverges since $\sum_{p\equiv1\ (3)}1/p=\infty$. Fix $m$ so that
this ratio exceeds $2j$. Removing every fiber of size less than $j$
therefore leaves $\gg P^3$ triples. Pairing each with a distinct
triple in its fiber gives $\gg P^3$ primitive points. The intervals
exclude every rational pairing plane except equality, which the
distinct triples avoid. If two different triples share a root,
deleting it leaves a nontrivial equal sum of two positive cubes.
Heath-Brown's bound \cite{c3:hbcubicsurface} gives
$O_\varepsilon(P^{4/3+\varepsilon})$ such ordered quadruples,
including their common multiples by summation. Hence only
$O_\varepsilon(P^{7/3+\varepsilon})$ pairs share a root.
Primitivity and positivity fix the projective scaling.
\end{proof}

The number of distinct targets with $j$ such representations is
$\gg P^{3-0.24871567}$ unconditionally. In this assertion and
both conditional bounds below, all roots can in addition be taken
$P^\eta$-smooth, for any prescribed $\eta>0$.

To impose smoothness take a smaller positive $\eta$ if necessary
for \cite[Theorem 1.2]{c3:wooley}. There is a positive-density
subset of the smooth numbers with the natural relative frequencies
in every fixed residue class, also in each prescribed interval.
For completeness choose finitely many disjoint prime intervals
$(P^{a_i},P^{(1+\omega)a_i})$ with $(1+\omega)a_i<\eta$,
$\sum a_i=1-\xi$, $0<\xi<\min a_i$, and
$(1+\omega)\sum a_i<1-\xi/2$. Choose their number large first,
then $\xi$ and the widths small. Let $d$ be a product of one
prime per interval and use the integers $dy\in PI$, with $y$
unrestricted. The cofactor is $O_I(P^\xi)$, smaller than every
marked prime, so each integer is counted once. For any fixed
modulus $q$ the marks are units and each root class has
$|I|P/(dq)+O(1)$ representatives per $d$. Mertens' theorem gives
$\sum_d1/d\to c_\eta>0$, and the total of the errors is
$O(P^{1-\xi/2})$. Thus the preceding first-moment sieve and the
choice of the auxiliary progression apply on this smooth subset;
only its fixed positive density changes.

Cauchy's inequality on the surviving triples now uses the mixed
sixth moment $O(P^{3+0.24871567})$ of
\cite[Theorem 1.2]{c3:wooley}. By orthogonality this bounds
the collision count when all three roots are smooth, since the
theorem allows one root on each side to be unrestricted.
The resulting target bound implies $\gg A^{0.91709477}$ distinct
targets up to $A$, by scaling the fixed intervals into $f\le A$.
The ordinary sixth moments under the Hooley--Heath-Brown and
Wang hypotheses give respectively $\gg_\varepsilon P^{3-\varepsilon}$
and $\gg P^3$ distinct targets. In each case one can delete
the $O_\varepsilon(P^{7/3+\varepsilon})$ targets having any two
distinct triples in these intervals sharing a root. In particular,
the positive-density hypotheses yield positive lower density for
every fixed multiplicity $j$, with all $3j$ smooth roots distinct.

\begin{corollary}\label{c3:growing-fibres}
There is $c>0$ such that $\gg A^{0.91709477}$ distinct targets up to $A$
each have at least
\[
                   \exp\{c\sqrt{\log A}/\log\log A\}
\]
representations in the prescribed intervals, with squarefree smooth
roots, pairwise coprime inside each triple and distinct across all
these representations.
\end{corollary}
\begin{proof}
Let $C:x^3+y^3+z^3=0$ in $\mathbb P^2$. This elliptic curve has
complex multiplication by $\mathbb Q(\sqrt{-3})$, as seen from
its order-three automorphism fixing a point. At a good prime $p$
put $a_p=p+1-\#C(\mathbb F_p)$. The CM equidistribution theorem
\cite[Corollary~3.8]{c3:fitesutherland} supplies $\gg y/\log y$
primes in $[y/2,y]$ with $a_p\le-\sqrt p$; one may apply the
theorem over the CM field and discard prime ideals of residue degree
greater than one. The affine zero fiber has $p^2-(p-1)a_p$ vectors.
Requiring all coordinates to be units removes $O(p)$ vectors,
since the three coordinate lines meet the projective curve in at
most nine points. Thus each such prime contributes a normalized
factor at least $1+1/(2\sqrt p)$.

Put $y=c_0\log P$, where $c_0>0$ is sufficiently small, and let
$m$ be the product of the selected primes. Write $G$ for the
product of these normalized factors. Then $m\le\exp(O(y))$ and
$G\gg\exp(c_1\sqrt y/\log y)$. Impose the target $0$ modulo $m$
and units for all three roots there. The smooth construction above
gives $\gg P^3G/m$ triples with the global squarefreeness and
coprimality conditions. We detail uniformity in this growing
modulus. The cofactor intervals in the construction have lengths
$\gg P^{\xi/2}$; each residue class therefore contains $\gg P/(dq)$
choices uniformly for moduli $q\le P^{\xi/4}$.
Impose the sieve up to a large fixed $L$ together with $m$, taking
$c_0$ small enough. For $L<p\le P^{1/8}$ outside $m$, count
unrestricted roots modulo $mp^2$. Every root class in a fixed
box contains $O(P/(mp^2))$ integers, since $mp^2=o(P)$.
The Chinese remainder theorem bounds the excluded mass by
$O((P^3G/m)\sum_{p>L}p^{-2})$. The remaining primes cost
$O(P^{3-1/8})$ by unrestricted packing. First choose $L$ large,
then $c_0$ small; the lower bound follows.

There are $O(P^3/m)$ target values. Deleting fibers of size
less than a sufficiently small fixed multiple of $G$ preserves
$\gg P^3G/m$ triples. Cauchy's inequality and the smooth moment
now give $\gg P^{3-0.24871567}/m^2$ such targets.
Since $3-0.24871567>3(0.91709477)$, $c_0$ can also absorb
the factor $m^2$ in this strict margin. Delete the
$O_\epsilon(P^{7/3+\epsilon})$ shared-root targets as above.
The result follows on taking $P$ a fixed multiple of $A^{1/3}$
small enough for the intervals.
\end{proof}

The size of this congruence gain has matching order in the
exponent. For every modulus $m$ and every target, the normalized
density of primitive root triples modulo $m$ is at most
$\exp(O(\sqrt{\log(3m)}/\log\log(9m)))$.
Indeed at a prime $p>3$ the primitive fibers lift smoothly.
The projective curve bound above gives $1+O(p^{-1/2})$ for the
zero fiber. For nonzero targets the cubic Gauss estimate gives
$1+O(1/p)$, or one may subtract the infinity curve from the
smooth cubic surface. The prime $3$ contributes a bounded factor
at arbitrary depth by fixing a unit variable and lifting its cube
on each unit class. Finally sum $p^{-1/2}$ separately for primes
up to $\log(3m)$ and for larger prime divisors of $m$, using the
prime number theorem and their product bound.

If common prime factors within individual triples are allowed, the
same number of distinct targets has at least
\begin{equation}\label{c3:imprimitive-many}
                    \exp\{c\log A/\log\log A\}
\end{equation}
representations with all roots still distinct, squarefree and
$P^\eta$-smooth. In the proof take instead every sufficiently large
prime $p\in[y/2,y]$ and prescribe target zero modulo $p^3$.
The solutions with all three coordinates units give normalized
mass $1+O(p^{-1/2})$, by smooth lifting from the plane cubic.
The disjoint solutions with all three coordinates of valuation
one give exactly $(1-1/p)^3$: every such triple has sum
divisible by $p^3$. Thus each prime contributes at least $3/2$
for large $p$. The product modulus remains $\exp(O(y))$, while
the gain is $\exp(\gg y/\log y)$. The same sieve can keep
pairwise coprimality at every unselected prime. The shared-root
bound already includes nonprimitive quadruples.

For unrestricted roots this scale also has matching order in the
exponent for single-congruence gains. At depth $j\ge3$ an all-divisible
triple contributes the normalized density at depth $j-3$ for the
scaled target when it is divisible by $p^3$. Primitive lifts
contribute at most $1+O(p^{-1/2})$ per recursion away from $3$,
and a bounded amount at $3$. The normalized count is therefore at
most a bounded factor times
$\tau(m)\prod_{p\mid m}(1+O(p^{-1/2}))$.
The elementary maximal-order bound for $\tau$ gives
$\exp(O(\log(3m)/\log\log(9m)))$.

\subsection{A range of the square-divisor sieve}\label{c3:sieve}
Let $F(\mathbf h)=h_1^3+\cdots+h_6^3$ and let $D=F^\vee$ be the primitive
equation of the dual hypersurface. For $1\le L\le Z^{3/2}$ define
\[
 B(Z,L)=\#\{\mathbf c\in[-Z,Z]^6\cap\mathbb Z^6:
           p^2\mid D(\mathbf c)\text{ for some prime }L\le p<2L\}.
\]
Wang's conjecture SF asks for $B(Z,L)\ll Z^6L^{-\delta}$ with fixed
$\delta>0$ throughout this range \cite[Conjecture 1.6]{c3:wang}.
The following estimate uses the original cubic to control divisibility of
its dual.
\begin{proposition}\label{c3:sieveprop}
For each $0<\eta<1/3$ there is $\delta(\eta)>0$ such that
\[
 B(Z,L)\ll_\eta Z^6 L^{-\delta(\eta)}
 \quad (1\le L\le Z^{4/3-\eta}).
\]
\end{proposition}
\begin{proof}
We may omit finitely many primes, adjusting the implied constant.
Call a residue $\mathbf c\bmod p$ good if $D(\mathbf c)=0$, all its
coordinates are nonzero and $\nabla D(\mathbf c)\ne0$ modulo $p$.
The other zeros lie on a fixed set of codimension at least two, since
$D$ is absolutely irreducible and none of the coordinates or the fixed
nonzero derivative $\partial_1D$ is divisible by $D$.
For $p\le 2Z$ a set of this kind has $O(p^4)$ residue vectors and contributes
$O(Z^6/p^2)$ integer vectors. For $L\ge Z/2$ the union of the bad sets
contributes $O(Z^5\log(2Z))$. To see this, change coordinates over
$\mathbb Z$ by a fixed unimodular matrix so that $D$, viewed as a polynomial
in the last variable, has nonzero constant leading coefficient. Such a
choice is possible by homogeneity, choosing a primitive last direction
where $D$ is nonzero. Take the resultant of $D$ and
$(\prod c_i)\partial_1D$ in that variable. The $O(Z^4)$ base points
annulling this nonzero resultant cost $O(Z^5)$. Elsewhere the prime
must divide its nonzero value of polynomially bounded size. There are
$O(\log(2Z))$ possible primes and $O(1+Z/p)$ last coordinates per prime.

At a good residue one partial derivative is invertible. On fixing the
other five coordinates the resulting simple roots give $O(1)$ classes
modulo $p^2$ by Hensel's lemma. In particular
\[
 B(Z,L)\ll Z^6/L+Z^5L\qquad(1\le L\le Z).
\]
This suffices for $L\le Z^{9/10}$. We treat the remaining range by Poisson
summation. Use $e_m(t)=\exp(2\pi i t/m)$ and write
\[
 T_p(\mathbf h)=\sum_{\substack{\mathbf c\bmod p^2\, :\,D(\mathbf c)=0\bmod p^2\\
                      \mathbf c\bmod p\ {\rm good}}}
                        e_{p^2}(\mathbf h\cdot\mathbf c).
\]
This transform is real by symmetry. If $p\nmid\mathbf h$, summation over
the lifts of a fixed good residue $\bar{\mathbf c}$ gives zero unless
$\bar{\mathbf h}$ is normal to the dual at $\bar{\mathbf c}$. In the
latter case it gives $p^5$ times the character at any one lift. The Gauss-map argument of \cite[proof of Proposition 9.11]{c3:wang} identifies both varieties on their smooth contact loci. The normal
line recovers the unique contact point $[\bar{\mathbf h}]$ on $F=0$,
and $[\bar{\mathbf c}]=[\nabla F(\bar{\mathbf h})]$.
For an elementary verification, away from the coordinate hyperplanes
adjoin square roots of the $c_i$. The factors of $D$ there are the sums
$\sum_i\pm c_i^{3/2}$, with sign patterns identified under simultaneous
negation \cite[Sect.~1]{c3:wang}. Smoothness means that precisely one factor
vanishes. Its gradient is proportional to the corresponding signed vector
of square roots. The normal and contact identities follow. Conversely,
for $F(\mathbf h)=0\bmod p$ the residue $\nabla F(\mathbf h)$ is good
except possibly on a fixed set $E$ of dimension at most four in the
six-dimensional affine space. One may take the part of $F=0$ where some
coordinate or some nonempty proper subsum of the $h_i^3$ vanishes. The
equations of the exceptional set are independent of $p$, and
$\#E(\mathbb F_p)\ll p^4$.

Lift a nonzero contact point to $\tilde{\mathbf h}\equiv\mathbf h\pmod p$
with $F(\tilde{\mathbf h})=0\pmod{p^2}$. Then
$\mathbf c_0=\nabla F(\tilde{\mathbf h})$ lifts a contact hyperplane
with $D(\mathbf c_0)=0\pmod{p^2}$, since the polynomial
$D(\nabla F)$ is divisible by $F$ away from finitely many primes.
Euler's identity and Taylor expansion give
$\mathbf h\cdot\mathbf c_0\equiv F(\mathbf h)\pmod{p^2}$.
The $p-1$ nonzero scalar multiples of the hyperplane modulo $p$
therefore supply the Ramanujan factor
$p^5(p\mathbf1_{p^2\mid F(\mathbf h)}-1)$.
If $p\nmid F(\mathbf h)$ there is no contact contribution. If the
contact hyperplane is excluded, the contribution is zero. Hence for
$p\nmid\mathbf h$ we have the one-sided bound
\begin{equation}\label{c3:ram}
 T_p(\mathbf h)\le
 p^4\sum_{a\bmod p^2}^{*}e_{p^2}(aF(\mathbf h))
  +p^5\mathbf1_{\bar{\mathbf h}\in E}.
\end{equation}
In particular, excluding a singular hyperplane costs only the size of
a negative Ramanujan term in this inequality.

For $\mathbf h=p\mathbf k$ and $p\nmid\mathbf k$, sum instead over the
scalar multiples of each good projective point of $D=0$. All but the
points on $\mathbf k\cdot\mathbf c=0$ give a factor $-1$, and those on
it give $p-1$. There are $O(p^4)$ points in total and $O(p^3)$ on any
such hyperplane, whence $|T_p(p\mathbf k)|\ll p^9$. For
$p^2\mid\mathbf h$ the trivial bound is $O(p^{10})$.
At either type of imprimitive frequency the Ramanujan term on the
right of \eqref{c3:ram} is nonnegative.

Majorize the original box by a constant times the Gaussian
$\exp(-\pi |\mathbf c|^2/Z^2)$. Poisson summation modulo $p^2$ gives
the factor $Z^6/p^{12}$ and the dual Gaussian
$\exp(-\pi |\mathbf h|^2/H^2)$, where $H=p^2/Z$. Its nonnegativity
permits the one-sided estimate \eqref{c3:ram}. Weighted residue counting
bounds the exceptional frequencies by $O(p^4(1+H/p)^6)$, the
$p$-divisible ones by $O((1+H/p)^6)$, and the $p^2$-divisible ones by
$O(1)$. After normalization all errors from these sets are
\[
 \ll (Z^6/p^{12})\{p^9(1+H/p)^6+p^{10}\}
 \ll Z^6/p^2+p^3 .
\]

The unrestricted Ramanujan term now factors into six one-dimensional sums.
Write $U(a)=\sum_n\exp(-\pi n^2/H^2)e_{p^2}(an^3)$. Its contribution
in absolute value is at most $(Z^6/p^8)\sum_a^*|U(a)|^6$.
Here $p\asymp L$, $Z^{9/10}\le L\le Z^{4/3}$ and $H\gg1$.
The classical cubic Weyl inequality, with the exact rational approximation
$a/p^2$, gives
\[
 \max_a^*|U(a)|^4\ll_\varepsilon Z^\varepsilon
        (H^4/p^2+H^3+p^2H).
\]
Partial summation and truncation at $H\log(2Z)$ justify the Gaussian form.
Also, by orthogonality,
\[p^{-2}\sum_{a\bmod p^2}|U(a)|^2\ll_\varepsilon Z^\varepsilon H,\]
since each unit $n$ has at most three cube-equivalent classes
modulo $p^2$, and the nonunits together contribute
$O((1+H\log(2Z)/p)^2)$. Substitution of $H=p^2/Z$ shows that the
Ramanujan contribution is $O_\varepsilon(Z^\varepsilon(Z^2p^2+Z^4))$.
Summing over the primes, and including the bad residues, bounds the
large-$L$ range by
\[
 B(Z,L)\ll_\varepsilon Z^\varepsilon
      (Z^6/L+Z^5+L^4+Z^2L^3+Z^4L).
\]
Each term saves a fixed power of $Z$ if $Z^{9/10}\le L\le Z^{4/3-\eta}$
and $\varepsilon$ is sufficiently small in terms of $\eta$. Together with
the slicing estimate this proves the proposition.
\end{proof}
The interval adjacent to $Z^{3/2}$ in SF remains open in this argument.
At the endpoint the unit-frequency sixth moment of the short cubic sum
with square denominator and the $L^4$ estimate for exceptional frequencies
both reach the ambient scale.

There is a further reduction under the Hasse--Weil hypotheses of
Hooley and Heath-Brown. Suppose the number $N_F(Y)$ of integral zeros
of $F$ in the signed box of radius $Y$ satisfies
$N_F(Y)\ll_\varepsilon Y^{3+\varepsilon}$, as in
\cite[Theorem 1.4]{c3:wang}. Then for every fixed $\eta>0$ the same
sieve saving holds up to $L\le Z^{3/2-\eta}$.

To verify this, keep the nonnegative Fourier majorant in the preceding
proof and sum the unit-frequency contribution over $p$ before applying
a mean value. Let $H=L^2/Z$ and
$U_s(\alpha)=\sum_{|n|\le s}e^{2\pi i\alpha n^3}$. The Gaussian sum
at scale $p^2/Z$ is an average of $U_s$ with the probability density
$2\pi s(Z/p^2)^2\exp(-\pi s^2(Z/p^2)^2)$ on $s>0$. For $L\le p<2L$
these densities are bounded by a constant multiple of the one at scale
$4H$. By convexity it suffices to apply the additive large sieve to
$U_s(\alpha)^3$ at the reduced fractions $a/p^2$. These fractions are
$\gg L^{-4}$ apart modulo one, and the polynomial has frequencies
in $[-3s^3,3s^3]$. The squared norm of its coefficient sequence is
$N_F(s)$, by changing signs on the last three variables. Integrating the
resulting large-sieve bound against the common Gaussian measure gives
\[
 \sum_{L\le p<2L}\sum_{a\bmod p^2}^*
       |U(a)|^6\ll_\varepsilon (L^4+H^3)H^{3+\varepsilon}
       \qquad(H\ge1).
\]
Here $U(a)$ on the left has its original $p$-dependent scale.
For $L\le Z^{3/2}$ one has $H^3\le L^4$. Multiplication by
$Z^6/L^8$ leaves $O_\varepsilon(Z^{3+\varepsilon}L^2)$ in the
range under consideration. The remaining errors are
$O_\varepsilon(Z^\varepsilon(Z^6/L+Z^5+L^4))$. These estimates save
a power up to $Z^{3/2-\eta}$. The constants deteriorate as $\eta$
decreases, and the endpoint $L\asymp Z^{3/2}$ still requires an estimate.

\subsection{Counting contact branches at the upper end of the sieve}
A different weight removes the exceptional primitive frequencies from
the Fourier analysis of Proposition~\ref{c3:sieveprop}. In particular
the top range of SF reduces to a sixth moment without an additional
error of order $L^4$.

\begin{proposition}\label{c3:contact}
For $Z^{9/10}\le L\le Z^{3/2}$ put
\[
 U_{p,Z}(a)=\sum_{n\in\mathbb Z}\exp(-\pi Z^2 n^2/p^4)e_{p^2}(an^3).
\]
Then
\[
 B(Z,L)\ll_\varepsilon Z^\varepsilon
       (Z^6/L+Z^5+L^{7/2}+ZL^3)
       +\frac{Z^6}{L^8}
         \sum_{L\le p<2L}\sum_{a\bmod p^2}^{*}|U_{p,Z}(a)|^6 .
\]
The four error terms save a fixed power of $Z$ throughout the stated
range. Thus a bound $\ll L^{8-\delta}$ for the displayed double sum,
uniformly from any fixed $Z^{3/2-\eta}$ to $Z^{3/2}$, would supply a
power saving on that part of SF.
\end{proposition}
\begin{proof}
Remove the singular and coordinate-zero residues in the original
$\mathbf c$-box by the resultant argument of
Proposition~\ref{c3:sieveprop}. For each remaining prime, majorize the
indicator of good solutions modulo $p^2$ by the number of solutions
in units $a,x_1,\ldots,x_5\bmod p^2$ of
\[
 x_1^3+\cdots+x_5^3+1=0,qquad
 \mathbf c=a(x_1^2,\ldots,x_5^2,1).
\]
Indeed $a=c_6$, and at each smooth point modulo $p$ the unique contact
direction is rational over the residue field. Normalize its last
coordinate to $1$. The square roots of $c_i/c_6$ with these fixed
residues lift uniquely. The local factorization of $D$ into signed sums
of the corresponding cubes shows that their chosen signed sum vanishes
modulo $p^2$, since the other factors are units. Thus the parametrization
counts each originally retained solution. It also counts the distinct
contact branches above some singular residues, with nonnegative weights.

Write $T'_p(\mathbf h)$ for the Fourier transform of this new weight.
For $p\nmid\mathbf h$, summing the lifts of a fixed residue of $(a,\mathbf x)$
requires both
$\sum_{i\le5}h_i x_i^2+h_6=0$ and
$h_i x_i=b x_i^2$ for $1\le i\le5$ and some scalar $b$ modulo $p$.
These conditions follow by varying $a$ and then the tangent vectors to
$\sum x_i^3+1=0$. Since the $x_i$ are units, the first condition gives
$h_6=b$ and the other five give $h_i=b x_i$. Every $h_i$ must thus
be a unit. The residue of $\mathbf x$ is unique when it exists, and
it exists precisely when $F(\mathbf h)=0\bmod p$. For each nonzero
scalar residue there are $p^5$ lifts of the parameters. As in
\eqref{c3:ram}, lifting the contact vector and using Euler's identity
shows that its scalar phase vanishes modulo $p^2$ precisely when
$p^2\mid F(\mathbf h)$. The sum over scalar residues therefore gives
the exact identity for primitive frequencies
\[
 T'_p(\mathbf h)=\mathbf1_{p\nmid\prod h_i}\,p^4
                     \sum_{a\bmod p^2}^{*}e_{p^2}(aF(\mathbf h)).
\]
Thus the six Gaussian Fourier sums factor with the restriction $p\nmid n$.
Dropping this restriction in the sixth moment costs at most a fixed
factor times the unrestricted sixth moment plus
$O(p^2(1+H/p)^6)$ at each prime, where $H=p^2/Z$. After multiplication
by $Z^6/p^8$ this latter term is $O(1+Z^6/p^6)$.

At $\mathbf h=p\mathbf k$, $p\nmid\mathbf k$, summation over scalars gives
\[
 T'_p(p\mathbf k)=p^5(pN_{\mathbf k}-N),
\]
where $N$ is the number of projective points of the cubic fourfold $F=0$
with all coordinates nonzero, over $\mathbb F_p$, and $N_{\mathbf k}$
also imposes $\sum k_i x_i^2=0$. The Lang--Weil estimates give
$N=p^4+O(p^{7/2})$. If
\[
                 \prod_{\varnothing\ne J\subset\{1,\ldots,6\}}
                           \sum_{j\in J} k_j^3\ne0\pmod p ,
\]
the cubic and quadric meet transversely, since at a singular intersection
their gradients would force the nonzero $x_i$ to be a common scalar
times $k_i$, and the cubic equation gives just such a vanishing subsum.
The complete intersection is smooth, geometrically connected in this
dimension, and hence geometrically irreducible. Its homogeneous ideal
contains no linear equation, so that removing the coordinate walls
costs $O(p^2)$ points. The uniform bounded-degree Lang--Weil estimate
gives $N_{\mathbf k}=p^3+O(p^{5/2})$, hence
$|T'_p(p\mathbf k)|\ll p^{17/2}$. For any nonzero $\mathbf k\bmod p$
the quadric cannot contain the cubic, and the dimension bound gives
$|T'_p(p\mathbf k)|\ll p^9$. Finally $|T'_p(0)|\ll p^{10}$.

The Gaussian on $\mathbf k$ has scale $H/p=p/Z$, smaller than $p$.
The fixed nonzero polynomial defining the exceptional $\mathbf k$
has $O((1+H/p)^5(\log(2Z))^{O(1)})$ zeros in the effective box,
by the elementary polynomial zero bound on a grid of distinct residues. The tail past $(1+H/p)\log(2Z)$ is negligible. The zero mode
and all $p^2$-divisible frequencies cost $O(p^{10})$ in the weighted
transform sum. After Poisson normalization the imprimitive contributions
at one prime are therefore
\[
 \ll_\varepsilon \frac{Z^{6+\varepsilon}}{p^{12}}
        \{p^{10}+p^{17/2}(1+H/p)^6+p^9(1+H/p)^5\}
 \ll_\varepsilon Z^\varepsilon(Z^6/p^2+p^{5/2}+Zp^2).
\]
Summation, together with the primitive-frequency identity and the
initial resultant estimate, proves the proposition.
\end{proof}

\subsection{An additive question at square modulus}
\label{c3:dual-energy}
The remaining short sum in the square-divisor problem has a useful
equivalent description at the upper end of the prime range. Let
$p>3$, $H=p^2/Z$, and put
\[
 U(a)=\sum_{p\nmid n}\exp(-\pi n^2/H^2)e_{p^2}(a n^3).
\]
Choose one $r$ from each of the two square classes in
$(\mathbb Z/p^2\mathbb Z)^\times$. For $v\in\mathbb Z$ prime to $p$,
consider every $y\bmod p^2$ with
$3r y^2\equiv-v\pmod{p^2}$, giving this pair weight
$\exp(-\pi v^2/Z^2)$. Define $E_3(Z;p^2,r)$ to be the sum of the
products of the six weights over the ordered sextuples of these pairs
for which
\[
             y_1^3+y_2^3+y_3^3\equiv
                    y_4^3+y_5^3+y_6^3\pmod{p^2}.
\]
The Fourier transform on the two square classes gives
\begin{equation}\label{eq:c3-square-energy}
 \frac{Z^6}{p^8}\sum_{a\bmod p^2}^{*}|U(a)|^6
             \ll\sum_r E_3(Z;p^2,r).
\end{equation}
Indeed write $a=rt^2$ and complete the Gaussian sum modulo $p^2$.
For a fixed $v$, summing the original variable over the $p$ lifts of
its residue class kills it unless $3ax^2+v=0\pmod p$. For unit
$x$ this condition forces $p\nmid v$. Each stationary residue has
a unique lift modulo $p^2$ solving the same equation. With $y=tx$,
the complete sum is therefore
\[
             p\sum_{3ry^2=-v\ (p^2)}
                    e_{p^2}(-2rt^{-1}y^3).
\]
Poisson summation contributes $H/p^2$ and the dual weight displayed
above. Summing the sixth power over the invertible $t$ and enlarging
to all additive frequencies gives at most $p^2 E_3(Z;p^2,r)$.
Since $H/p=p/Z$, the factors cancel in \eqref{eq:c3-square-energy}.

At $L\asymp Z^{3/2}$, an estimate
$\sum_{p\asymp L}\sum_r E_3(Z;p^2,r)\ll Z^{6-\eta}$ for some
$\eta>0$ would give the required power saving at that scale by the
contact parametrization above. It asks for a saving over $Z^{9/2}$
in the energy per prime on average, up to logarithmic factors.
This formulation retains the original tangencies: a tuple satisfying
the cubic relation gives a normal hyperplane by taking the six
squares. The gain sought concerns how often such square classes lie
in short intervals. The additive energies of modular roots studied
in \cite{c3:roots} give a comparison in prime modulus. Here both
the cube on the root variables and the second $p$-adic digit enter
the congruence.

Wang \cite[Proposition~9.11]{c3:wang} showed that the square-divisor
conjecture for the quaternary Fermat dual implies it for all diagonal
cubics in at least four variables. Here the fourth moment also identifies
the scales and the contribution lost on reduction modulo $p$.
Give the residue $y^3$ the sum of the Gaussian weights
of its pairs $(v,y)$ above, and write $E_j(q,r)$ for the $j$th additive
energy after projecting these weighted residues to $\mathbb Z/q$, where
$q=p$ or $p^2$. The mass for each $r$ is $O(Z)$, since each admissible
$v$ has two roots. Young's convolution inequality implies
\[
 E_3(p^2,r)\ll Z^2 E_2(p^2,r).
\]
Thus $\sum_{p\asymp L,r}E_2(p^2,r)\ll Z^{4-\eta}$ for some $\eta>0$
would suffice at $L\asymp Z^{3/2}$. At least one of the two square
classes has mass $\gg Z$, since every $v$ prime to $p$ occurs for one
of the classes. Cauchy's inequality on the $p$ residue classes gives
\[
          \sum_{L\le p<2L}\sum_r E_j(p,r)
               \gg_j Z^{2j}/\log L\qquad(j=2,3).
\]
Consequently even an optimal estimate that retains only the first
$p$-adic digit cannot give the required power saving by this test.

The part of $E_2(p^2,r)$ supported on pair sums divisible by $p$ is
$O(Z^2)$. Indeed, for fixed $y_1$ the condition
$y_1^3+y_2^3=0\bmod p$ permits at most three residues for $y_2$.
Each fixes $v_2\bmod p$. The total Gaussian weight of its lifts $v_2$
is bounded since $Z<p$, and the root lifting that residue $y_2$ is
unique for each such $v_2$. The total mass of the relevant pairs is
$O(Z)$, whose square bounds this part of the energy. After summation
over $p$ it is already $O(LZ^2)$.

For a unit pair sum $t$, put $b=-3r$ and $\lambda=b^{-3}$ modulo $p^2$.
The mass at $t$ is exactly the weighted number of pairs of integers
$v_1,v_2$ prime to $p$ on the congruence
\begin{equation}\label{eq:cube-unit-pair-curve}
 (t^2-\lambda(v_1^3+v_2^3))^2
                   =4\lambda^2 v_1^3v_2^3\pmod{p^2}.
\end{equation}
In fact, with $X_i=y_i^3$ we have $X_i^2=\lambda v_i^3$.
Conversely, set
$X_1=(t^2+\lambda v_1^3-\lambda v_2^3)/(2t)$ and $X_2=t-X_1$.
The congruence gives both required squares. Then $y_i=bX_i/v_i$
recovers a unique pair with $b y_i^2=v_i$ and $y_i^3=X_i$.
Thus the remaining fourth moment concerns the distribution of the
short-interval counts on the curves \eqref{eq:cube-unit-pair-curve},
as $t$ varies through units at precision $p^2$.

For illustration, take the sharp cutoff $|v|\le\lfloor p^{2/3}\rfloor$
and $r=1$. The counts can be reproduced by lifting the first digit of
a square root and convolving integer multiplicities on $\mathbb Z/p^2$.
\[
 \begin{array}{c|r|r|r|r}
       p & \lfloor p^{2/3}\rfloor &\text{root pairs}&E_2(p^2,1)&E_3(p^2,1)\\\hline
       499&62&124&47148&43547740\\
       743&82&164&81708&100896080\\
       1009&100&220&147036&269787280
 \end{array}
\]
The calculation covers individual primes and uniform interval weights.
A proof of SF by this route requires an averaged estimate uniform at
the upper end of the prime range.

\begin{proposition}\label{c3:independent-radicals}
At $L\asymp Z^{3/2}$ the terms of
$\sum_{p\asymp L}\sum_r E_2(p^2,r)$ for which the four nonzero
integers $v_i$ are linearly dependent in
$\mathbb Q^\times/\mathbb Q^{\times2}$ contribute
$O_\varepsilon(Z^{7/2+\varepsilon})$. The part in which the
four-variable dual form vanishes over $\mathbb Z$ satisfies the
same estimate. Consequently only independent squareclasses remain
in the fourth-moment criterion.
\end{proposition}
\begin{proof}
Truncate the Gaussians at $W=Z\log(2Z)$, at an error
$O_A(Z^{-A})$ for any fixed $A$. The omitted mass is bounded
by summing the four weights without the congruence, since each
unit $v_i$ admits at most two roots for fixed $p,r$. Define the
degree-twelve integer polynomial
\[
 D_4(\mathbf v)=\prod_{\epsilon_2,\epsilon_3,\epsilon_4=\pm1}
    \big(\sqrt{v_1^3}+\epsilon_2\sqrt{v_2^3}
              +\epsilon_3\sqrt{v_3^3}+\epsilon_4\sqrt{v_4^3}\big).
\]
Independent sign changes leave this product invariant. Thus it
is a polynomial in the $v_i$ and is independent of the branches.
Every energy tuple has $p^2\mid D_4(\mathbf v)$: substitute
$v_i=-3r y_i^2$ and use homogeneity and the cubic relation.

If $D_4(\mathbf v)\ne0$, there are $O(1)$ possible dyadic-range
primes dividing its value, since $|D_4(\mathbf v)|\ll W^{12}$.
Each permits $O(1)$ choices of all the branches. Dependence gives
a nonempty subset of $j$ coordinates whose product is a rational
square. Its absolute value is an integer square $m^2\le W^j$.
For fixed $m$ there are $O_\varepsilon(W^\varepsilon)$ ordered
factorizations into $j$ signed integers, by the divisor bound.
The total for this subset is
$O_\varepsilon(W^{4-j/2+\varepsilon})$.

If $D_4(\mathbf v)=0$, write $v_i=d_i n_i^2$ with signed
squarefree $d_i$ and $n_i\ge1$. Square roots from distinct rational
squareclasses are linearly independent, as is seen by applying
the character projections of their multiquadratic extension over
$\mathbb Q$. The vanishing signed sum of the $v_i^{3/2}$ therefore
splits into blocks with equal $d_i$, each of size at least two.
Two blocks of size two force equality of the two $n_i$ within each
block and allow $O(W^2)$ vectors. A single block requires
$\sum_i\epsilon_i n_i^3=0$. For $1\le n_i\le T$ this has
$O_\varepsilon(T^{2+\varepsilon})$ solutions for each sign pattern.
With both signs present one can place a pair of the same sign on
one side and fix the two remaining variables; a positive integer
$m$ has $O_\varepsilon(m^\varepsilon)$ expressions as $a^3+b^3$,
since $a+b$ divides $m$ and its value fixes the quadratic equation
for $a,b$. The all-one-sign pattern contributes nothing. Summing
$T=(W/|d|)^{1/2}$ over the common signed squarefree $d$ gives
$O_\varepsilon(W^{1+\varepsilon})$. Thus the integer-zero part has
$O(W^2)$ vectors, each contributing $O(L)$ over primes and branches.
Both parts satisfy the asserted bound.
\end{proof}

\subsection{Nearby factorable denominators in the square-divisor sieve}
\label{c3:nearby-denominators}
At the upper end of the sieve, a reduced denominator $p^2$ lies in
the classical Weyl range for the cubic sum. A second rational
approximation can have a denominator with three coprime factors,
to which Heath-Brown's estimate applies \cite[Theorem~2]{c3:hb-factor}.
The following result isolates the frequencies where this approach
can still encounter the Weyl exponent.

Put $Q=p^2$, $u=10^{-6}$ and $M=\lfloor Q^{1/1000}\rfloor$.
Let $\mathcal D(Q)$ consist of the integers $q=q_1q_2q_3$ with
\[
 Q^{2/7}<q_1\le2Q^{2/7},\qquad Q^{5/21}<q_2\le2Q^{5/21},\qquad
 Q^{10/21-u}<q_3\le2Q^{10/21-u},
\]
where $q_1,q_2$ are primes, $q_3$ is squarefree and the factors are
pairwise coprime. Call a unit $a\bmod Q$ exceptional if the interval
$[1,M]$ misses all residues $aq\bmod Q$ for $q\in\mathcal D(Q)$.

\begin{proposition}\label{c3:factorable-numerators}
For all sufficiently large $p$ there are $O(Q^{1-10^{-5}})$
exceptional numerators. If $p\le 2Z^{3/2}$ and $a$ is
nonexceptional, the unit Gaussian cubic sum $U(a)$ in
Section~\ref{c3:dual-energy} satisfies
$|U(a)|\ll p^{1-u/4}$, uniformly for $H=Q/Z\ge1$.
For $Z^{7/5}\le L\le Z^{3/2}$ one consequently has
\[
 \sum_{L\le p<2L}\ \sum_{\substack{a\bmod p^2\ {\rm unit}\\
                             a\ {\rm nonexceptional}}}|U(a)|^6
                    \ll_\varepsilon L^{8-u/2+\varepsilon}.
\]
\end{proposition}
\begin{proof}
Give each denominator weight $w(q)$ equal to its number of
factorizations of the displayed type. All these denominators are
less than $Q$ when $Q$ is sufficiently large. Write
$A=\sum w(q)=Q^{1-u+o(1)}$. The prime number theorem and the
elementary squarefree density estimate give the lower bound, since
discarding multiples of $q_1$ or $q_2$ from the third interval costs
$O(q_3/q_2+q_3/q_1+1)$ per prime pair, with $q_3$ here denoting its
scale. Also $\sum w(q)^2\ll_\varepsilon A Q^\varepsilon$ by the
divisor bound. Regard $w$ as a function on the group $G$ of units
modulo $Q$. For each nonprincipal character $\chi$ of $G$,
\begin{equation}\label{eq:c3-factorable-character}
 |\widehat w(\chi)|\ll A Q^{-1/30}.
\end{equation}
For a primitive character use the $r=2$ Burgess estimate
$\ll_\varepsilon B^{1/2}Q^{3/16+\varepsilon}$ on intervals of
length $B$ \cite{c3:burgess}. An imprimitive nonprincipal character
comes from a character modulo $p$, to which the same estimate applies
with the smaller modulus. Write the squarefree indicator as
$\sum_{d^2\mid n}\mu(d)$. With $Y=Q^{10/21-u}$, the terms
$d\le Q^{1/25}$ contribute
$\ll_\varepsilon Y^{1/2} Q^{3/16+\varepsilon}$, and
$d>Q^{1/25}$ costs $O(YQ^{-1/25}+\sqrt Y)$. The removed multiples
of either prime are covered by the preceding elementary bound.
Summing over the two prime coordinates proves
\eqref{eq:c3-factorable-character}.

Let $J$ be the indicator of $[1,M]$ in $G$, and put
$R(a)=\sum_{q}w(q)J(aq)$. Its mean on $G$ is
$\mu=AM/|G|$. Parseval gives
\[
 \frac{1}{|G|}\sum_{a\in G}\frac{|R(a)-\mu|^2}{\mu^2}
  =\frac{1}{A^2M^2}\sum_{\chi\ne1}
                         |\widehat w(\chi)|^2|\widehat J(\chi)|^2 .
\]
The characters with $|\widehat J(\chi)|\le M^{99/100}$ contribute
$O_\varepsilon(Q^{u+\varepsilon}M^{-1/50})$. There are
$O_\varepsilon(Q^{0.022+\varepsilon})$ remaining characters.
Indeed, the products of $999$ integers at most $M$ lie below $Q$,
so orthogonality and the divisor bound give
$\sum_\chi|\widehat J(\chi)|^{1998}\ll_\varepsilon
 Q^{1+\varepsilon}M^{999}$. Markov's inequality gives the stated
count. Their normalized contribution is
$O_\varepsilon(Q^{0.022-2/30+\varepsilon})$ by
\eqref{eq:c3-factorable-character}. Taking $\varepsilon>0$
sufficiently small proves the first assertion.

For a nonexceptional $a$, write $aq-bQ=h$, $1\le h\le M$.
Reducing $b/q$ removes $g=(b,q)\mid h$. Both primes $q_1,q_2$
exceed $M$, so the reduced denominator has the same two prime
factors and third factor $q_3/g$, still squarefree. Denote it by
$q'=q/g$, and write $q'_3=q_3/g$. First consider the sharp sum
$S(a/Q,t)=\sum_{n\le t}e_Q(an^3)$ at lengths
$Q^{2/3-u/3}\le t\le C Q^{2/3}\log Q$, where $C$ is fixed.
Then $t\le q'\le t^{3/2}$ and
$t\le\min(q_1q'_3,q_2q'_3)$. Heath-Brown's theorem gives
\[
 |S(a/Q,t)|\ll_\varepsilon
 \left(1+\frac{t^3h}{Qq}\right) Q^\varepsilon
  \big(t^{1/2}q_1^{1/2}+t^{1/4}(q'q_2)^{1/4}
                   +t^{1/4}(q')^{1/4}(q'_3)^{1/8}\big).
\]
The three terms in parentheses are
$\ll_\varepsilon Q^{10/21+\varepsilon}$ and the approximation
factor is $\ll_\varepsilon Q^{1/1000+u+\varepsilon}$, hence the
whole expression is $O(Q^{0.48})$.
At smaller $t$ the ordinary cubic Weyl inequality for the exact
fraction $a/Q$ gives
\[
 |S(a/Q,t)|\ll_\varepsilon Q^\varepsilon
       \big(tQ^{-1/4}+t^{3/4}+t^{1/4}Q^{1/4}\big)
       \ll Q^{1/2-u/8}.
\]
Here $H\le 2^{2/3}Q^{2/3}$, so partial summation applies these
bounds to the Gaussian. Removing multiples of $p$ costs
$O(1+H/p)$. This proves the asserted bound.

The fourth-power version of \eqref{eq:c3-square-energy} gives
\[
 \sum_{L\le p<2L}\sum_{a\bmod p^2}^*|U(a)|^4
 \ll (L/Z)^4 L^2 \sum_{p\asymp L}\sum_r E_2(Z;p^2,r)
 \ll_\varepsilon L^{6+\varepsilon}.
\]
For the last estimate, away from $D_4=0$ each quadruple admits
only $O(1)$ primes in this range, since $|D_4|\ll
(Z\log(2Z))^{12}$ after truncation. On $D_4=0$ the count
$O_\varepsilon(LZ^{2+\varepsilon})$ from
Proposition~\ref{c3:independent-radicals} suffices. Multiplying by
the squared pointwise bound proves the sixth-moment estimate.
\end{proof}

The unresolved contribution in Proposition~\ref{c3:contact} can
therefore be restricted to numerators for which every one of the
residues $h a^{-1}$, $1\le h\le M$, avoids the indicated three-factor
set. The size estimate for that exceptional set concerns counting
measure. The remaining sieve question concerns its fourth-power
weight.

\subsection{The prime terms of the reciprocal-coefficient mean}
Write $V$ for the Fermat cubic fourfold, $D=F^\vee$, and let
$\mu_{\mathbf c}(n)$ be the coefficients of $1/L(s,V_{\mathbf c})$,
with the middle cohomology normalized to weight zero. The prime terms
in Wang's smoothed mean square admit the following unconditional bound.
\begin{proposition}\label{c3:primemean}
For $2\le M\le Z$ and arbitrary coefficients $|b_p|\le1$,
\[
 \sum_{\substack{\|\mathbf c\|_\infty\le Z\\D(\mathbf c)\ne0}}
 \left|\sum_{M\le p<2M} b_p\mu_{\mathbf c}(p)\right|^2
       \ll Z^6 M.
\]
The same estimate holds for the prime coefficients of
$\Phi_{\mathbf c,1}$ in \cite[Conjectures~6.4, 6.7]{c3:wang}.
\end{proposition}
\begin{proof}
Omit the finitely many primes of bad reduction of $V$. On
$\mathbb F_p^6$ put
\[
 f_p(\mathbf c)=p^{-3/2}\sum_{x\in V(\mathbb F_p)}
             (\mathbf1_{\mathbf c\cdot x=0}-1/p).
\]
For $p\nmid D(\mathbf c)$ the trace formula for the smooth threefold
gives $\mu_{\mathbf c}(p)=f_p(\mathbf c)+O(p^{-1/2})$:
the projective-space contribution is
$\#\mathbb P^3(\mathbb F_p)$ and
$\# V(\mathbb F_p)=\#\mathbb P^4(\mathbb F_p)+O(p^2)$.
The normalized prime coefficient of $\Phi_{\mathbf c,1}$ differs
from $\mu_{\mathbf c}(p)$ by $O(p^{-1/2})$ as well. These errors
have sum $O(\sqrt M)$ for each $\mathbf c$.

We first bound the mean square with $f_p$ on the whole box. The
additive Fourier expansion is exact:
\[
 f_p(\mathbf c)=p^{-5/2}
   \sum_{\substack{\mathbf h\bmod p,\ \mathbf h\ne0\\ F(\mathbf h)=0}}
              e_p(\mathbf c\cdot\mathbf h).
\]
There are $O(p^5)$ frequencies. Parseval, followed by covering the
box by complete residue boxes, gives a diagonal contribution
$O(Z^6 M)$, since $p\le2Z$. Majorize the box by a Gaussian of scale
$Z$ in treating the off-diagonal. For distinct $p,q$ the Fourier
terms are parameterized by integer vectors
$\mathbf t\equiv q\mathbf h-p\mathbf k\pmod{pq}$. Such a vector
specifies at most one pair of frequencies and must satisfy
$pq\mid F(\mathbf t)$. The zero vector cannot occur. By Poisson
summation the absolute off-diagonal contribution is at most
\[
 \ll Z^6M^{-5}
 \sum_{p\ne q}\sum_{\substack{\mathbf t\ne0\\pq\mid F(\mathbf t)}}
             \exp(-c\|\mathbf t\|^2/T^2),\qquad T=M^2/Z .
\]
For $F(\mathbf t)\ne0$ and $\|\mathbf t\|\le 2^j T$, at most
$O((1+j)^2)$ ordered prime pairs divide $F(\mathbf t)$, as $T\le M$.
The resulting sum, with the prefactor, is $O(Z^6M^{-5}T^6)$, using
Gaussian decay also when $T<1$. This is $O(M^7)$. At $F(\mathbf t)=0$
Hua's eighth moment and the fourth moment give
$O_\varepsilon((2^jT)^{7/2+\varepsilon})$ nonzero integral vectors
in this ball whenever it has radius at least one. The total is
$O_\varepsilon(Z^6M^{-3}T^{7/2+\varepsilon})$ for $T\ge1$, with
an arbitrarily high power of $T$ available for $T<1$. Thus these
terms also contribute $O(Z^6M)$.

It remains to compare at bad primes of the section. For
$\mathbf c\not\equiv0\pmod p$, every hyperplane section has at worst
isolated singularities, since the Gauss map
$[x_i]\mapsto[x_i^2]$ has finite fibres. The point-count bound for
singular complete intersections \cite{c3:hooleysingular} gives
$f_p(\mathbf c)\ll\sqrt p$ there. At $\mathbf c\equiv0\pmod p$
one has $f_p(0)\ll p^{5/2}$. The actual normalized local
coefficients are bounded in every case by \cite[(3.2)--(3.3)]{c3:wang}.
The correction at a bad prime is therefore bounded by
\[
 C\big(\sqrt M\,\mathbf1_{p\mid D(\mathbf c)}
              +M^{5/2}\mathbf1_{p\mid\mathbf c}\big).
\]
If $2M\le\sqrt Z$, complete residue counting at $p$ and at $pq$
for distinct primes gives
$\#\{\mathbf c:p\mid D(\mathbf c)\}\ll Z^6/p$ and the two-prime
bound $\ll Z^6/(pq)$. Here $D$ has $O(p^5)$ zeros modulo $p$.
For the coordinate-divisibility indicators use $p^6$ and
$(pq)^6$ in the same bounds. Expanding the square, their total is
$O(Z^6M)$. In the complementary range each $\mathbf c$ with
$0<|D(\mathbf c)|\ll Z^{\deg D}$ has only $O(1)$ prime divisors
in $[M,2M)$. Cauchy's inequality then reduces the correction to the
same one-prime counts, still valid for $p\le2Z$.
\end{proof}
The full reciprocal mean extends to $M=Z^3$ and includes composite
indices. In the off-diagonal calculation the unrestricted geometric
trace equals $f_p(0)\asymp p^{5/2}$ at the origin. With positive
prime weights this single excluded vector already contributes a
quantity of order $M^7/(\log M)^2$ to the extended mean. Estimates
beyond the box scale must remove this concentration before taking
absolute Fourier values.

\subsection{Rough moduli in the reciprocal mean}
The prime mean extends to all moduli with large prime factors,
including prime powers. Write $P^-(1)=\infty$ and $P^-(n)$ for
the least prime factor of $n>1$. Let $a_{\mathbf c}(n)$ denote
either $\mu_{\mathbf c}(n)$ or the coefficient of
$\Phi_{\mathbf c,1}$ in \cite[(1.6)]{c3:wang}.
\begin{theorem}\label{c3:roughmean}
For each fixed $0<\eta<1$, $2\le M\le Z$, and $|b_n|\le1$,
\[
 \sum_{\substack{\|\mathbf c\|_\infty\le Z\\D(\mathbf c)\ne0}}
 \left|\sum_{\substack{M\le n<2M\\ P^-(n)\ge Z^\eta}}
         b_n a_{\mathbf c}(n)\right|^2 \ll_\eta Z^6 M .
\]
\end{theorem}
\begin{proof}
We may take $Z$ sufficiently large in terms of $\eta$ that the
rough primes have good reduction for $V$. For squarefree $n$ put
$f_n=\prod_{p\mid n} f_p$, with $f_p$ as in the proof of
Proposition~\ref{c3:primemean}. First consider integer affine
progressions $\mathbf c=\mathbf a+d\mathbf z$. Uniformly in
$\mathbf a,d$, $1\le B\le X\le 2Z$, and $|u_n|\le1$, we claim
\begin{equation}\label{c3:roughgeom}
 \sum_{\|\mathbf z\|_\infty\le X}
  \left|\sum_{\substack{B\le n<2B,\ n\text{ squarefree}\\
                 (n,d)=1,\ P^-(n)\ge Z^\eta}}
        u_n f_n(\mathbf a+d\mathbf z)\right|^2
             \ll_\eta X^6 B.
\end{equation}
The claim permits boxes of side comparable with $X$ and intervals
with lower endpoint between $1/2$ and $1$. For each $p$ the mean
of $|f_p|^2$ is $p^{-5}\#\{\mathbf h\ne0:F(\mathbf h)=0\}
=1+O(p^{-2})$ by the smooth projective count. It bounds all
normalized Fourier coefficients of $|f_p|^2$ in absolute value.
Thus the product of these bounds over any set of primes is bounded.
The diagonal in \eqref{c3:roughgeom} is $O(X^6B)$ by completing
residue boxes modulo $n\le2X$.

For distinct $n,m$ write $g=(n,m)$, $n=gu$, $m=gv$, and
$Q=B/g$, so $u,v\asymp Q$. At primes of $g$ take the preceding
bound for the Fourier coefficients of $|f_p|^2$. At each prime
of $uv$ the normalized Fourier coefficient is $p^{-5/2}$ on the
nonzero zeros of $F$ and zero elsewhere. The translation in the
progression introduces phases of absolute value one and its
dilation permutes these frequencies by scalars. Poisson summation
against a Gaussian majorant of the box now gives an off-diagonal
contribution bounded by
\[
 X^6\sum_{g} Q^{-5}
   \sum_{u,v}\ \sum_{\substack{\mathbf t\in\mathbb Z^6\setminus\{0\}\\
                                  uv\mid F(\mathbf t)}}
             \exp(-c\|\mathbf t\|^2/T^2),
 \qquad T=B^2/(gX).
\]
Indeed the period is $guv$, and each $\mathbf t$ specifies all
local frequencies uniquely. At least one prime is exclusive to
one of $n,m$, so $\mathbf t=0$ makes no contribution. Here $u,v$
are coprime and all divisors displayed are $Z^\eta$-rough. When
$0<|F(\mathbf t)|\ll (2^jT)^3$ the number of eligible ordered pairs
$(u,v)$ is at most $\exp(O_\eta(1+j))$: assign each rough prime
dividing $F(\mathbf t)$ to $u$, $v$, or neither. The Gaussian absorbs
this bound. The nonzero-value contribution is at most
\[
  \ll_\eta X^6\sum_{g}Q^{-5}T^6
  \ll_\eta B^7\sum_{\substack{g\le2B\\P^-(g)\ge Z^\eta}}1/g
  \ll_\eta B^7 .
\]
For the last sum, the Euler product is bounded by Mertens' theorem
on $[Z^\eta,4Z]$. The same Gaussian bound for lattice points works
for $T<1$, by the exponential decay away from the origin.

When $F(\mathbf t)=0$, the sixth moment bounds the Gaussian count
by $O(T^{7/2+1/4})$ for $T\ge1$, and by an arbitrarily high power
of $T$ below one. There are $O(Q^2)$ pairs $u,v$ per $g$. Since
$T\le Q$ and $Q\ge1/2$, the contribution is
$\ll X^6\sum_{g\le2B}(B/g)^{3/4}\ll X^6B$.
This proves \eqref{c3:roughgeom}.

At good primes of the section the trace formula gives
$a_{\mathbf c}(p)=f_p(\mathbf c)+h_p$ with $h_p$ independent of
$\mathbf c$ and $h_p=O(p^{-1/2})$. This holds in both choices
of coefficient, since their difference at $p$ is independent
of $\mathbf c$. The local weight bound in \cite[(3.2)--(3.3)]{c3:wang}
and the singular point-count bound used above yield
\[
 a_{\mathbf c}(p)=f_p(\mathbf c)+h_p+e_p(\mathbf c),\qquad
 |e_p(\mathbf c)|\ll \sqrt p\,\mathbf1_{p\mid D(\mathbf c)}
                       +p^{5/2}\mathbf1_{p\mid\mathbf c}.
\]
Set $e_1=1$ and $e_s=\prod_{p\mid s}e_p$ for squarefree $s$.
At $D(\mathbf c)\ne0$ there are $O_\eta(1)$ rough squarefree
$s$ with $e_s(\mathbf c)\ne0$, since
$|D(\mathbf c)|\ll Z^{\deg D}$. First expand the squarefree
products using $f_p+e_p$. Cauchy's inequality in this bounded
set of $s$ reduces the desired square mean to
\[
 \ll_\eta\sum_{s}\sum_{D(\mathbf c)\ne0}|e_s(\mathbf c)|^2
   \left|\sum_{r\asymp M/s} u_{sr}f_r(\mathbf c)\right|^2 ,
\]
where $r$ is squarefree, rough, and coprime to $s$; $s$ is squarefree
and rough. On each residue class modulo $s$ bound $|e_s|^2$ by the
product of the preceding indicator majorants squared. Their sum over
the residue classes is at most $C^{\omega(s)}s^6$ by the Chinese
remainder theorem: for a prime $p$ the discriminant has $O(p^5)$
zeros, and the all-zero class has weight $O(p^5)$. Applying
\eqref{c3:roughgeom} in each class with $X\asymp Z/s$ therefore costs
$\ll_\eta Z^6 (M/s) C^{\omega(s)}$ per $s$. The sum of
$C^{\omega(s)}/s$ over these rough $s\le2M$ is bounded in terms
of $\eta$. Endpoint scales of order one are covered by the
one-term bound in \eqref{c3:roughgeom}.

The constants $h_p$ may now be restored by expanding over their
squarefree product $h_t$. Minkowski's inequality costs a factor
$\sum_t|h_t|/\sqrt t\ll_\eta1$ in the square root of the mean,
by the same bounded Euler product. This establishes the result for
squarefree moduli, with arbitrary further support restrictions on
the coefficients. Finally factor each $n=br$ uniquely with $b$
squarefull, $r$ squarefree and $(b,r)=1$. For $n\le2Z$ rough,
$|a_{\mathbf c}(b)|\ll_\eta1$ by the uniform local weight bound
and the bounded degrees of the Euler factors. A final application
of Minkowski's inequality costs
\[
 \sum_{\substack{b\le2M\ \text{squarefull}\\P^-(b)\ge Z^\eta}}
 b^{-1/2}\le
 \prod_{Z^\eta\le p\le 2M}
 \left(1+\frac{p^{-1}}{1-p^{-1/2}}\right)\ll_\eta1 .
\]
This completes the proof.
\end{proof}

\begin{theorem}\label{c3:fullshortmean}
There is an absolute constant $C$ such that for $2\le M\le Z$ and
$|b_n|\le1$ the full, unrestricted short mean satisfies
\[
 \sum_{\substack{\|\mathbf c\|_\infty\le Z\\D(\mathbf c)\ne0}}
 \left|\sum_{M\le n<2M} b_n a_{\mathbf c}(n)\right|^2
          \ll Z^6M\big(\log(2Z)\big)^C .
\]
\end{theorem}
\begin{proof}
We give the two estimates that replace roughness in the preceding
proof. First, for each fixed $A$ there is a $C_A$ such that
\begin{equation}\label{c3:divcubic}
 \sum_{\substack{\|\mathbf t\|_\infty\le T\\F(\mathbf t)\ne0}}
   \tau(|F(\mathbf t)|)^A \ll_A T^6\big(\log(2T)\big)^{C_A}
       \qquad(T\ge1).
\end{equation}
For $1\le k\le 6T^3$ one has
$\tau(k)^A\ll_A\sum_{d\mid k,\, d\le T}\tau(d)^{B_A}$ with
a constant $B_A$. Indeed the prime factors exceeding $T$ contribute
a bounded factor when $T\ge 2$. Group the remaining prime factors,
with multiplicity, successively into maximal blocks of product at
most $T$. The product of two consecutive blocks exceeds $T$, so
the number of blocks is bounded. Submultiplicativity of $\tau$
proves the assertion. The number of zeros of $F$ modulo $d$ is
$O(d^5)$. At $p>3$ the nonzero zeros lift smoothly, with
$p^5+O(p^3)$ primitive solutions at level one. If all coordinates
are divisible by $p$, scaling out $p$ reduces the level by three.
In the ratio of the zero count to $p^{5j}$ this term has multiplier
$p^{-3}$ for $j\ge3$, giving a local bound $1+O(p^{-2})$ at every
level. At $p=3$, fix one unit coordinate modulo three. Its cube
on $x=r+3z$ is $r^3+9(r^2z+3rz^2+3z^3)$, with the expression in
parentheses a permutation at every level. This bounds the primitive
local ratio uniformly; the same scaling bounds the full ratio.
Completing boxes modulo $d\le T$ and summing $\tau(d)^{B_A}/d$
by its Euler product proves \eqref{c3:divcubic}.

In \eqref{c3:roughgeom} the number of possible pairs $u,v$ at
$F(\mathbf t)\ne0$ is bounded by $\tau(|F(\mathbf t)|)^2$.
Using \eqref{c3:divcubic} and $\sum_{g\le2B}1/g\ll\log(2B)$
therefore gives the same estimate with a fixed power of $\log(2Z)$
on the right and without the lower bound on the prime factors.
The argument at $F(\mathbf t)=0$ already sums over every $g$.

Secondly, the local error of the preceding proof satisfies
\begin{equation}\label{c3:localerr}
  \sum_{\mathbf a\bmod p}\ 
     \sup_{\substack{\mathbf c\equiv\mathbf a\ (p)\\D(\mathbf c)\ne0}}
           |e_p(\mathbf c)|^2 \ll p^5 .
\end{equation}
Outside a closed subset of codimension at least two in
$\mathbb A^6$, a point of the dual gives a section with one ordinary
double point and no other singularities. To see this, remove the
singular locus of the dual and its intersections with the coordinate
hyperplanes. Gauss duality gives a unique contact there. If its
coordinates are all nonzero, the diagonal Hessian of $F$ is
nondegenerate, and the radical on its tangent hyperplane is the
line of contact, proving ordinarity. These statements hold in all
sufficiently large characteristics. The contact point is rational
over $\mathbb F_p$ by uniqueness. Projection from this point
writes the section as $x_0 Q_2(\mathbf x)+Q_3(\mathbf x)=0$ with
$Q_2$ a smooth quadric in $\mathbb P^3$. The curve $C=(Q_2=Q_3=0)$
is smooth: a dependence between its two gradients would give a
second singular point of the original section. As a complete
intersection it is geometrically connected, hence geometrically
irreducible, and has genus four. Counting points on the projection
lines gives precisely
\[
 1+\#\mathbb P^3(\mathbb F_p)-\#Q_2(\mathbb F_p)
       +p\#C(\mathbb F_p)=\#\mathbb P^3(\mathbb F_p)+O(p^{3/2}).
\]
Consequently $e_p=O(1)$ on the ordinary one-node locus, which has
$O(p^5)$ points. On the remaining $O(p^4)$ nonzero residue vectors
of the dual one has $e_p=O(\sqrt p)$ by the bound for isolated
singularities. The origin contributes $O(p^5)$ to the square.
This proves \eqref{c3:localerr}.

For a fixed nonsingular section apply weighted Cauchy to the divisor
sum in the preceding proof, with weights $1/s$:
$\sum_{s\mid D(\mathbf c)}1/s\ll\log(2Z)$. The complementary
weight $s$ is paid for by \eqref{c3:localerr}: the squared
majorants over residue classes modulo $s$ now have sum at most
$C^{\omega(s)}s^5$ instead of $C^{\omega(s)}s^6$. The remaining
sum $\sum_{s\le2M}C^{\omega(s)}/s$ contributes a power of the
logarithm. Restoring the deterministic prime terms and the
squarefull parts costs additional fixed log powers via the same
Euler products, since all prime-power coefficients of
$a_{\mathbf c}$ are uniformly bounded: the reciprocal local
factors are polynomials of fixed degree with bounded reciprocal
roots. The finitely many excluded primes may likewise be separated
at bounded cost. This proves the theorem.
\end{proof}

\subsection{A log-free mean below the box scale}
The two corrections in Wang's $\Phi_{\mathbf c,1}$ control the
contribution from prime squares and the mean at primes. Write
$A_{\mathbf c}(n)$ for its coefficients.
\begin{theorem}\label{c3:logfreeall}
Fix $0<\sigma<1$. For arbitrary $|b_n|\le1$ and
$2\le M\le Z\exp(-(\log(2Z))^\sigma)$ one has
\[
 \sum_{\substack{\|\mathbf c\|_\infty\le Z\\D(\mathbf c)\ne0}}
   \left|\sum_{M\le n<2M} b_n A_{\mathbf c}(n)\right|^2
              \ll_\sigma Z^6 M .
\]
\end{theorem}
\begin{proof}
For all sufficiently large primes put $G_{p,0}=1$, $G_{p,1}=f_p$
and, for $e\ge2$, let $G_{p,e}$ equal $A_{\mathbf c}(p^e)$ on
the good-reduction residue classes and zero elsewhere. These
functions depend only on $\mathbf c\bmod p$. There is a fixed $R$
such that $G_{p,e}$ and $A_{\mathbf c}(p^e)$ vanish for $e>R$;
those with $e\ge2$ are uniformly bounded. Write
$G_n=\prod_{p^e\Vert n}G_{p,e}$. On taking expectations over residue
vectors one has
\begin{equation}\label{c3:momentseuler}
 \mathbb E G_{p,1}=0,\qquad
 \mathbb E |G_{p,1}|^2=1+O(p^{-2}),\qquad
 \mathbb E G_{p,2}=O(p^{-1}).
\end{equation}
Only the last assertion requires an additional input.
Equations (6.5)--(6.10) of \cite{c3:wang} are unconditional
finite-field estimates. The unrestricted reciprocal has mean
$1+O(p^{-1})$ at $p^2$, also on restricting to good reduction
by zero extension. The factor $1/\zeta(2s)$ subtracts the
good-reduction probability $1+O(p^{-1})$. The shift of $1/L(s,V)$
contributes $O(p^{-1})$, since its prime coefficient is
$O(p^{-1/2})$ and the restricted prime mean of the unrestricted
reciprocal is $O(p^{-1/2})$.

We claim that $G_n$ satisfies a log-free geometric mean on a
progression $\mathbf c=\mathbf a+d\mathbf z$ with
$\|\mathbf z\|_\infty\le X$, for $n\asymp B$ coprime to $d$
and $B\ge1/2$, provided that $1\le X\ll Z$ and
$X/B\gg\exp(\tfrac13(\log(2Z))^\sigma)$. Its bound is
$O_\sigma(X^6B)$, uniformly in the progression and in
bounded coefficients.

For the zero Fourier mode put
$m_p(e,f)=|\mathbb E G_{p,e}\overline{G_{p,f}}|$.
Set $h_p(e,f)=m_p(e,f)$ except for
$h_p(1,1)=\max(0,m_p(1,1)-1)$. Then
\[
 m_p(e,f)\le\sum_{j=0}^{\min(e,f)}h_p(e-j,f-j),\qquad
 \prod_p\sum_{e,f\ge0}\frac{h_p(e,f)}{p^{(e+f)/2}}<\infty.
\]
Indeed \eqref{c3:momentseuler} treats $e+f\le2$ and higher
terms are bounded and have total degree at least three. With
$h(u,v)=\prod_p h_p(v_p(u),v_p(v))$, the double sum of the
complete absolute means over $n,m\asymp B$ is at most
\[
 \sum_{u,v}h(u,v)\#\{j:ju,jv\asymp B\}
    \ll B\sum_{u,v}\frac{h(u,v)}{\sqrt{uv}}\ll B .
\]
This gives $O(X^6B)$ after Gaussian Poisson summation.

For a nonzero Fourier mode let $U$ be the product of primes of
total multiplicity one in $nm$, and $V$ the product of the other
primes dividing it. Then $UV^2\le4B^2$, the period is $UV$, and
the normalized Fourier coefficient at the frequency specified
by $\mathbf t$ has modulus at most $U^{-5/2}C^{\omega(V)}$.
In fact the primes of $U$ require $U\mid F(\mathbf t)$ by the
Fourier formula for $f_p$; each prime of $V$ is bounded by
Cauchy's inequality and the bounded second moments. The translation
and the unit dilations in the progression preserve these facts.
For a fixed $V$ there are at most $C^{\omega(V)}$ choices of its
pairs of exponents. Take $U\asymp U_0$ dyadically with
$U_0\ll B^2/V^2$ and put $T=U_0 V/X$.
When $F(\mathbf t)\ne0$, the number of choices for $U$ and its
allocation between $n,m$ is $\ll\tau(|F(\mathbf t)|)^2$.
By \eqref{c3:divcubic} the resulting sum is
\[
 \ll X^6(\log(2Z))^C\sum_{V\le2B} C^{\omega(V)}
         \sum_{U_0} U_0^{-5/2} T^6
 \ll B^7(\log(2Z))^{C'} .
\]
Gaussian decay handles all dyadic dilates, also when $T<1$.
At $F(\mathbf t)=0$, $\mathbf t\ne0$, use the cubic sixth
moment $O(T^{7/2+1/4})$ for $T\ge1$, and Gaussian decay below
one. The unrestricted count of allocations of $U\asymp U_0$
is $O(U_0\log(2U_0))$. Writing $\alpha=1/4$, the resulting bound
is
\[
 \ll X^6(\log(2Z))^C X^{-7/2-\alpha} B^{4+2\alpha}
   \sum_{V\le2B} C^{\omega(V)}V^{-1/2-\alpha}
 \ll X^6B(B/X)^{7/2+\alpha}(\log(2Z))^{C'} .
\]
The gap between $B$ and $X$ absorbs the logarithms in both
estimates, proving the geometric mean.

At a prime of good reduction the precise difference
$A_{\mathbf c}(p)-f_p(\mathbf c)$ is $p^{-5/2}$:
the fixed fourfold term is canceled by the shifted reciprocal,
and $\#\mathbb P^4/p-\#\mathbb P^3=1/p$. For every $e\ge1$
write
\[
 A_{\mathbf c}(p^e)=G_{p,e}+\mathbf1_{e=1}p^{-5/2}
                      +E_{p,e}(\mathbf c).
\]
The last term is supported on $p\mid D(\mathbf c)$.
Its squared supremum over lifts has residue mean $O(1/p)$
by \eqref{c3:localerr} for $e=1$ and by boundedness for $e\ge2$.
The same bound holds for
$\mathbb E(\mathbf1_{p\mid D}|G_{p,e}|^2)$ for every $e$,
by the three strata in the proof of \eqref{c3:localerr}.
The deterministic constants may be set aside: expanding over
the product $t$ of their primes and using Minkowski in the final
estimate costs at most $\sum_t t^{-3}<\infty$.

Expand the remaining polynomial over the part $d$ that uses
$E_{p,e}$. Its complementary factor $r$ is coprime to $d$.
Put $\theta=1-\sigma/4$ and
$Y=\exp(\tfrac12(\log(2Z))^\sigma)$. On $D(\mathbf c)\ne0$,
weighted Cauchy in $d$ costs
\[
  \sum_{d\mid |D(\mathbf c)|^R} d^{-\theta}
  \ll_\sigma 1+
  \sum_{d'\le Y}(d')^{-\theta}
           \mathbf1_{\operatorname{rad}(d')\mid D(\mathbf c)} .
\]
Indeed the tail above $Y$ is at most
$Y^{-\sigma/4}\prod_{p\mid D}(1-p^{-1+\sigma/2})^{-1}$.
Since $|D|\ll Z^{\deg D}$, splitting at $p=\log(2Z)$ bounds
the logarithm of this product by $O_\sigma((\log(2Z))^{\sigma/2})$.
The complementary Cauchy weight on the square is $d^\theta$.

For each $d'$ in the last sum, separate from $r$ the factor $k$
supported on primes of $d'$. Its exponents are at most $R$, so
Cauchy over $k$ costs at most $(R+1)^{\omega(d')}$. In each class
modulo $l=\operatorname{rad}([d,d'])$ the remaining $G$-polynomial
has coprime indices and length $B\asymp M/(dk)$. Its box scale
is $X\asymp Z/l$, with $l\le 2MY<Z$ for large $Z$, hence $X/B\gg Z/(MY)$ as required. Residues
with no possible index may be omitted. The squared local
majorants from $E_d$, $G_k$, and the indicated discriminant
divisibilities have joint mean at most $C^{\omega(l)}/l$ in
these classes: at each prime of $d$ use the error bound, and
at each other prime of $d'$ use the restricted $G$ bound.
Dropping $1/k$ in the resulting estimate leaves a total at most
\[
 \ll_\sigma Z^6 M
  \sum_{d,d'\ge1}
   \frac{C^{\omega(d)+\omega(d')}}{
        d^{1-\theta}(d')^\theta\operatorname{rad}([d,d'])}.
\]
This Euler product converges: the lowest nonconstant terms at $p$
are $O(p^{-2+\theta}+p^{-1-\theta}+p^{-2})$. The omitted fixed
primes and their powers contribute bounded factors. The theorem
follows.
\end{proof}

For the Mellin formulation one can allow all the test functions of
\cite[Conjecture 6.4]{c3:wang}, with its stated norm, in the range
$N\le Z^{5/6-\varepsilon}$. Indeed Cauchy's inequality in the
inversion integral gives a weight at $n$ bounded by a constant
times the square root of that norm times $\min(1,(N/n)^2)$, by
using the lines $0$ and $2$. Below
$Z\exp(-(\log(2Z))^{1/2})$ sum Theorem~\ref{c3:logfreeall}
dyadically, obtaining $O(Z^3\sqrt N)$ in the square root of the
mean, with the norm factor suppressed. Above that scale a dyadic
block at $B$ has square-root mean
\[
 \ll_\delta \min(Z^{1/2}B^{3+\delta}, Z^3 B^{1+\delta}).
\]
For the first bound enlarge the family box to $Y=B^{1+\delta/4}$.
The coefficients depend only on the projective hyperplane. On
each primitive integer ray of height $h\le Z$ the ratio of the
two multiplicities is $\lfloor Z/h\rfloor/\lfloor Y/h\rfloor
\ll Z/Y$, so Theorem~\ref{c3:logfreeall} gains this factor in
the square mean. The other bound follows from the uniform divisor
estimate on $A_{\mathbf c}(n)$. Using the first up to $B=Z^{5/4}$
and the second thereafter costs
$O_\delta(N^2 Z^{7/4+2\delta})$. Choosing $\delta$ small in
terms of $\varepsilon$ proves the stated range.

Under automorphy and GRH for the hyperplane-section $L$-functions
this last all-test-function conclusion holds for $N\le
Z^{1-\varepsilon}$. To bound the tail, insert a fixed smooth dyadic
partition $\psi(n/B)$. Write the inverse Mellin weight on the line
$\Re s=2$ and the dyadic cutoff on a second line of integration
$z$. Shift the latter to $\Re(s+z)=1/2+\delta$. The standard
right-half-plane bound for the reciprocal in
\cite[Proposition 3.2]{c3:wang}, with the fixed shifted factors
absolutely convergent there, gives a bound
$O_\delta(Z^{O(\delta)} B^{1/2+\delta}(N/B)^2)$ per section,
times the square root of the test norm. Indeed the Mellin transform
of $\psi$ is rapidly decreasing and
$\int|f(2+it)|(1+|t|)^\delta\,dt$ is controlled by that norm.
Use the unconditional estimate on dyadic blocks up to
$Z^{1-\varepsilon/2}$; beyond this point $B/N\ge
Z^{\varepsilon/2}$ absorbs the losses.

\subsection{Prime coefficients beyond the box scale}
The large geometric trace at the affine origin can be removed before
estimating the Fourier covariance. A trilinear exponential sum then
extends the prime range through the box scale.
\begin{theorem}\label{c3:longprimemean}
For every $\varepsilon>0$ the conclusion of
Proposition~\ref{c3:primemean} holds throughout
$2\le M\le Z^{7/6-\varepsilon}$. If the sixth moment for six
cubes is $O_\delta(H^{3+\delta})$ for each $\delta>0$, it holds
throughout $2\le M\le Z^{3/2-\varepsilon}$. The latter hypothesis
follows from the Hasse--Weil assumptions of
\cite[Theorem 1.4]{c3:wang} for this family.
\end{theorem}
\begin{proof}
Take $Z\le M$, since the shorter range is known. Choose a
nonnegative Schwartz majorant $W$ of the unit box whose Fourier
transform is smooth and compactly supported. For instance square
the modulus of an inverse Fourier transform of a smooth bump and
dilate. In estimating the geometric mean we sum $W(\mathbf c/Z)$
only over $\mathbf c\ne0$. The diagonal is $O(Z^6M)$.
At $p\nmid D(\mathbf c)$ use the bounded smooth trace. For each
$0<|D(\mathbf c)|\ll(2^jZ)^{\deg D}$ only $O(1+j)$ primes
$p\asymp M$ divide $D(\mathbf c)$, and their nonzero residue
vectors give $|f_p|\ll\sqrt p$. The vectors on $D=0$ and
$\|\mathbf c\|\le B$ number $O(B^3(\log(2B))^3)$: group their
nonzero coordinates by signed rational squareclass. In a vanishing
signed sum of the radicals $c_i^{3/2}$ each group has size at least
two by character independence. A group of size $j\ge2$ gives
$O(B^{j/2}\log(2B))$ choices by summing over the squarefree part.
Thus the nonzero residues of these vectors cost
$O(Z^3(\log(2Z))^3 M^2)$. Finally, $\mathbf c\equiv0\pmod p$,
$\mathbf c\ne0$, contributes
$O_A(M^6(Z/M)^A)$ for any $A>6$, by Schwartz decay.
These estimates also justify replacing the geometric sum by the
actual reciprocal on the desired nonsingular vectors:
the deterministic discrepancy is $O(\sqrt M)$ per vector and
the bad-prime discrepancies away from zero residues have squared
sum $O(M)$ per vector. Zero residues in the original box are
handled by the same bound, with a bounded number of bad primes
per nonsingular vector.

Put $T=M^2/Z$ and $C=T^3/M^2$. For distinct primes $p,q\asymp M$,
Poisson summation on all $\mathbf c$ transforms the covariance to
\[
 Z^6(pq)^{-5/2}
 \sum_{\substack{\mathbf t\in\mathbb Z^6,\ pq\mid F(\mathbf t)\\
          \mathbf t\not\equiv0\ (p),\ \mathbf t\not\equiv0\ (q)}}
                \widehat W(Z\mathbf t/(pq)).
\]
Deleting the two restrictions on the zero residues costs at most
$O(M^2(1+T/M)^6)$ in the unnormalized sum over pairs, since the
weight has compact support. This is $O(M^6)$ in our ranges.
The full congruence is now $F(\mathbf t)=pqm$. Insert a fixed
smooth compactly supported weight $w(m/C)$ which is identically
one on all possible values; $\|\mathbf t\|\ll T$.
Write the weighted count as an integral over the unit circle.

For the minor arcs this integral separates the cubic sum in
$\mathbf t$ from the trilinear sum over $p,q,m$.
To justify the separation, write $v=\log(pq/M^2)$ and expand
the smooth function $e^{-5v/2}\widehat W(e^{-v}\mathbf t/T)$
in a Fourier series in $v$, after inserting a compact cutoff
identically one on $[0,\log 4]$ and extending periodically.
Every resulting weight in $\mathbf t/T$ is smooth on a fixed
compact set with all fixed derivative norms decaying faster than
any power of the Fourier index. The additional prime weights
have bounded modulus. For a fixed resulting cubic weight
$w_1(\mathbf t/T)$ write
\[
 f_1(\alpha)=\sum_{\mathbf t}w_1(\mathbf t/T)e(\alpha F(\mathbf t)).
\]
Multidimensional partial summation and H\"older give
$\int|f_1(\alpha)|\,d\alpha\ll_\delta T^{\rho+\delta}$ with
$\rho=7/2$ unconditionally and $\rho=3$ under the stated sixth
moment. The constants depend on fixed derivative norms of $w_1$.
For arbitrary bounded prime-pair weights the trilinear sum obeys
\[
 |g(\alpha)|\ll_\delta T^\delta\,(M^2/r+1)(C+r\log(2r))
 \qquad (|\alpha-a/r|\le r^{-2},\ (a,r)=1).
\]
Indeed sum $m$ first and bound its transform by
$O(\min(C,\|\alpha pq\|^{-1}))$. There are at most two
ordered prime pairs of each product, and summation over all positive
integers $\ll M^2$ gives the displayed rational approximation bound.

Let $\Delta=C T^{-\eta}$, with $\eta>0$ small in terms of
$\varepsilon$. The major arcs have $r\le\Delta$ and
$|r\alpha-a|\le\Delta/T^3$. They are disjoint for large $Z$,
since $C\le T^{3/2}$ when $M\le Z^{3/2}$.
On their complement Dirichlet approximation gives
$\Delta<r\ll T^3/\Delta$. As $C\le M^2$ the preceding bound
yields $|g|\ll_\delta M^2 T^{\eta+\delta}$.
Thus the unnormalized minor integral over prime pairs is
$O_\delta(M^2 T^{\rho+\eta+2\delta})$.

On a major arc put $\beta=\alpha-a/r$. If $r\nmid pq$,
$\|\alpha pq\|\gg 1/r$, because
$|\beta|pq\ll\Delta/(Cr)$. Poisson summation in $m$ now gives
$O_A(C(1+C/r)^{-A})$ for every $A$. All such terms are negligible,
as $C/r\ge T^\eta$. The only possible noncentral resonant
denominators are $r=p$ and $r=q$; indeed $\Delta<pq$. At $r=p$
the $m$-sum is $C\widehat w(C\beta pq)$ up to a negligible error.
Completion of $f_1(a/p+\beta)$ modulo $p$, using the bound
$2\sqrt p$ for every linearly twisted complete cubic, gives
\[
 |f_1(a/p+\beta)|
     \ll T^6 p^{-3}(1+T^3|\beta|)^B .
\]
Here $B$ is fixed: the smooth transform decays after integration
by parts, $T\gg p$, and the derivative losses are polynomial in
$1+T^3|\beta|$. Integration against $C\widehat w(C\beta pq)$
costs $O(1/(pq))$ with the same polynomial. Summing $a,p,q$,
and likewise at $r=q$, gives $O(T^6/M^2)$.

For the central arc work with each original pair. The $m$-sum is
again $C\widehat w(C\beta pq)$ up to a negligible error. Its
rapid decay allows $|\beta|T^3\le T^{\eta/2}$, with a negligible
tail after rounding the exponent down if necessary. Poisson in
the six $\mathbf t$-variables on this smaller interval replaces
their sum by the corresponding real integral with a negligible
error, since its nonzero frequencies are nonstationary on the
support. Extending the $\beta$-integral to the real line and
applying Fourier inversion leaves exactly
\[
 \frac1{pq}\int_{\mathbb R^6}\widehat W(Z\mathbf t/(pq))\,d\mathbf t
           = (pq)^5 Z^{-6}W(0);
\]
the auxiliary $w(F(\mathbf t)/(Cpq))$ equals one on the support.
After the Poisson prefactor this is $(pq)^{5/2}W(0)$. The omitted
affine origin in the original covariance has weight
$W(0)f_p(0)f_q(0)$ with $f_p(0)=p^{5/2}(1+O(p^{-2}))$.
The total discrepancy after summing pairs is only $O(M^5)$.

The noncentral terms fit the target whenever
$T^{\rho+\eta+2\delta}\ll M^4$ and $T^6\ll M^8$,
since their prefactor is $O(Z^6M^{-5})$. With
$T=M^2/Z$ these hold throughout both claimed ranges, by choosing
$\eta,\delta$ sufficiently small. This proves the result.
\end{proof}

The upper end of the unconditional range can be made
$M\le Z^{7/6}(\log(2Z))^{2/3-\epsilon}$.
Vaughan's sixth-moment bound
\cite[Theorem~1.2]{c3:vaughan-trescubi} is
$O_\epsilon(T^{7/2}(\log(2T))^{-5/2+\epsilon})$.
It applies to the smooth six-variable weight by partial summation
and H\"older: the sixth power of a sum on any signed subinterval
is bounded in mean by the full collision count on a positive
interval of length $O(T)$, up to an absolute factor, on separating
the two signs. For an interval with a smaller endpoint scale,
monotonicity of the positive collision count also suffices.

In the proof above take $w(m/C)=\exp(-\pi (m/C)^2)$,
multiplying the smooth $\mathbf t$-weight for each pair by
$\exp(\pi(F(\mathbf t)/(Cpq))^2)$. This compensator and its
derivatives are bounded uniformly on the original compact
support after scaling by $T$, and have the same smooth Mellin
separation in $pq/M^2$. Put $L=K\sqrt{\log(2T)}$ and
$\Delta=C/L$, with $K$ a sufficiently large absolute constant.
The Gaussian gives an $O(T^{-A})$ bound with any prescribed fixed
$A$ at all nonresonant major centers after choosing $K$.
On minor arcs, when $|\alpha-a/r|\le r^{-2}$, split the product
indices into blocks of length at most $r/2$. Their frequencies
are separated by $1/(2r)$; the peak sum of the periodic
Gaussian in each block is $O(C+r)$. Thus
$|g|\ll (M^2/r+1)(C+r)\ll M^2 L$.
The minor integral becomes
$O_\epsilon(M^2T^{7/2}(\log(2T))^{-2+\epsilon})$.
At the central and noncentral resonant centers the Gaussian
transform has rapid decay, and Fourier inversion cancels the
compensator exactly at the central term. The other bounds
therefore remain those in the proof. The required minor
inequality now reads
$M^3\ll Z^{7/2}(\log(2Z))^{2-\epsilon}$, proving the refinement.

The exponent $3/2$ in the conditional range also describes a
limit of this use of absolute moments. Even away from the rational
pairing planes, the integer frequency $F(\mathbf t)=0$ supplies
$\gg T^3$ solutions in suitable fixed dilates, by the primitive
triple count above. The target inequality
$M^2T^3\ll M^6$ itself then requires $M\ll Z^{3/2}$.
An extension beyond this threshold must use cancellation in the
covariance or subtract an additional term before taking its
absolute value. The Fourier description gives one precise
subtraction. For each rational vector pairing space
$\Lambda\subset\{F=0\}$ of dimension three, the inverse transform
of $p^{-5/2}\mathbf1_{\Lambda\bmod p}$ is
$\sqrt p\,\mathbf1_{\Lambda^\perp\bmod p}$.
It vanishes on smooth hyperplane sections: a cubic threefold
containing a plane has singular points where the two normal
quadrics on the plane meet. Thus removing these finitely many
plane transforms, with their values at the zero frequency adjusted
separately, changes only the bad prime errors in the comparison
to the reciprocal coefficients. The nonplanar frequency zeros
are still present after this subtraction.

\begin{corollary}\label{c3:longrough}
Fix $0<\eta\le1$ and $0<\delta<\eta/6$. The rough-modulus
estimate of Theorem~\ref{c3:roughmean} holds up to $M\le
Z^{1+\delta}$, for both reciprocal coefficient families.
\end{corollary}
\begin{proof}
Write $Y=Z^\eta$ and take $M\ge Z$. A nontrivial squarefull part
$b$ has $b\ge Y^2$ and leaves length $M/b\le Z$; the bounded Euler
product in the proof of Theorem~\ref{c3:roughmean} disposes of all
these parts. Removing the deterministic prime discrepancies likewise
costs a bounded Euler product and leaves length $M/t\le Z$ when
$t>1$. For the remaining squarefree comparison, any nontrivial
product $d$ of bad-prime errors is at least $Y$. Use Cauchy over
the $O_\eta(1)$ possibilities at $D(\mathbf c)\ne0$ and
\eqref{c3:localerr}; in each class modulo $d$ enlarge the geometric
box to side $M/d$ for a polynomial at length $M/d$.
The total squared error is
\[
 \ll_\eta M^7\sum_{d\ge Y} C^{\omega(d)}/d^2
       \ll_\epsilon M^7 Y^{-1+\epsilon}\ll Z^6 M.
\]
It suffices to treat the squarefree geometric polynomials on the
whole punctured box with the weight of Theorem~\ref{c3:longprimemean}.

Pairs sharing a prime are bounded by M\"obius inversion of the
coprimality condition and the sum, over rough squarefree $d>1$,
of the means with $d\mid n,m$. Factor out $|f_d|^2$. It is
bounded by $C^{\omega(d)}$ times the product over $p\mid d$ of
$1+\mathbf1_{p\mid D}|f_p|^2$. For the empty subset of the
expanded product the shorter geometric mean at length $M/d$ in
the box of side $Z$ contributes $O_\eta(Z^6M/d)$.
For any nonempty subset with product $s$, the local factors have
joint residue mean $O(C^{\omega(s)}/s)$ by the three strata of
\eqref{c3:localerr}. Enlarge the box on each class modulo $s$
to side $M/s\ge M/d$. Summing the costs over $d=sh$ gives
$O_\epsilon(M^7Y^{-1+\epsilon})$ again, since the rough
harmonic sum over $h$ is bounded. The same estimates with a
Schwartz majorant follow by dyadic dilation.

For coprime $n,m$, use the proof of Theorem~\ref{c3:longprimemean}
with $pq$ replaced by $nm$. The multiplicity of each rough product
is $O_\eta(1)$, so the minor bound is unchanged. Dropping the
restrictions $\mathbf t\not\equiv0\pmod p$ for $p\mid nm$ costs
$O_\epsilon(T^6Y^{-5+\epsilon}+M^2(1+T^{7/2+\epsilon}))$
before the Poisson prefactor, with $T=M^2/Z$. Indeed at
$F(\mathbf t)\ne0$ all the other rough factors of $nm$ divide
this value in only $O_\eta(1)$ ways, and the nonzero multiples
of $p$ in the box number $O((T/p)^6)$ for $p\ll T$.
At $F=0$ there are $O(M^2/p)$ pairs containing $p$ and
$O_\epsilon(1+(T/p)^{7/2+\epsilon})$ corresponding multiples,
by scaling and the sixth moment.

An additional resonance can occur at any squarefree $r\mid nm$,
$r>1$, $r\le\Delta$. Completion of the cubic sum at these
centers costs $O_\eta(T^6 r^{-3}(1+r/T)^6)$ times a fixed
polynomial in the scaled arc coordinate. The smooth $m$-variable
of the trilinear argument (the free integer, distinct from the
modulus index here) integrates this polynomial at cost
$O(1/(nm))$. Since each pair has $O_\eta(1)$ possible divisors,
and $r\ge Y$, the aggregate resonant contribution is
\[
 \ll_\eta T^6\big(Y^{-2}+T^6/M^8\big).
\]
For $r>T$ this uses $r^4/T^6\le\Delta^4/T^6\le T^6/M^8$.
The central term cancels the affine origin with remaining error
$O_\eta(M^7Y^{-2})$ after normalization, by the product of
$1+O(p^{-2})$ over the rough factors. All unnormalized errors
above are $O(M^6)$ when $\delta<\eta/6$, using
$T^{7/2+\epsilon}\ll M^4$, $T^6Y^{-2}\ll M^6$ and
$T^{12}\ll M^{14}$. This proves the estimate.
\end{proof}

\end{document}